\documentclass[a4paper,12pt,reqno]{amsart}

\usepackage{amsmath,amssymb,amsthm}
\usepackage{color}
\usepackage{graphicx}
\usepackage{mathrsfs}
\usepackage{enumerate}
\usepackage{mathtools}
\usepackage[abbrev]{amsrefs}
\usepackage[T1]{fontenc}
\usepackage[colorlinks,
           linkcolor=blue,      
           anchorcolor=blue,  
           citecolor=red,       
            ]{hyperref}
\usepackage{amsthm}
\usepackage{amsmath,amssymb,amsthm}
\usepackage{latexsym}
\usepackage{color}
\usepackage{graphicx}
\usepackage{mathrsfs}
\usepackage{enumerate}
\usepackage[abbrev]{amsrefs}
\usepackage[T1]{fontenc}
\usepackage{mathtools}
\mathtoolsset{showonlyrefs=true}
\usepackage{appendix}

\allowdisplaybreaks[4]

\theoremstyle{plain}
\newtheorem{thm}{Theorem}[section]
\newtheorem{lemm}[thm]{Lemma}

\theoremstyle{definition}

\newtheorem{rem}[thm]{Remark}

\makeatletter

\@addtoreset{equation}{section}
\makeatother

\renewcommand{\div}{\operatorname{div}}

\renewcommand{\leq}{\leqslant}
\renewcommand{\geq}{\geqslant}
\newcommand{\ep}{\varepsilon}

\newcommand{\n}[1]{{\left\|#1\right\|}}

\newcommand{\f}{\frac}

\renewcommand{\sp}[1]{\left(#1\right)}

\newcommand{\R}{\mathbb{R}}
\newcommand{\Grad}{\nabla}

\newcommand{\p}{\partial}

\newcommand{\vu}{\mathbf{u}_{\varepsilon}}
\newcommand{\vv}{\mathbf{v}_{\varepsilon}}
\newcommand{\Om}{\Omega}

\newcommand{\Z}{\mathcal{Z}}

\begin{document}
\title[Low Mach number limit for two-fluid model]
{
Low Mach Number Limit and Convergence Rates for a Compressible Two-Fluid Model with Algebraic Pressure Closure
}

\author[Y.~Li]{Yang Li}
\author[M.~Luk\'a\v{c}ov\'a-Medvi\v{d}ov\'a]{M\'{a}ria Luk\'{a}\v{c}ov\'{a}-Medvi\v{d}ov\'{a}}
\author[E.~Zatorska]{Ewelina Zatorska}
\address[Y.~Li]{School of Mathematical Sciences, Anhui University, Hefei, 230601, People's Republic of China}
\email{lynjum@163.com}
\address[M.~Luk\'a\v{c}ov\'a-Medvi\v{d}ov\'a]{Institute of Mathematics, Johannes Gutenberg-University Mainz, Germany\\
RMU Co-Affiliate of Technical University Darmstadt}
\email{lukacova@uni-mainz.de}
\address[E.~Zatorska]{Mathematics Institute, University of Warwick, Zeeman Building, Coventry CV4 7AL, United Kingdom}
\email{ewelina.zatorska@warwick.ac.uk}
\date{\today} 

\begin{abstract}
We study the low Mach number limit for a viscous compressible two-fluid model with algebraic pressure closure in the three-dimensional torus $\mathbb{T}^3$. The pressure is determined implicitly through the densities of the two phases, which makes the singular limit substantially more delicate than for models with explicit pressure laws. Working in the framework of local-in-time strong solutions, we prove that, for well-prepared initial data, solutions to the rescaled compressible two-fluid system exist on a time interval independent of the Mach number and converge to the solution of the incompressible Navier--Stokes equations as the Mach number tends to zero. In addition, we establish explicit convergence rates for the densities and the velocity field. The proof relies on uniform high-order energy estimates and a relative energy argument adapted to the implicit structure of the pressure law. These results provide a rigorous justification of the low Mach number limit for the compressible two-fluid model with algebraic pressure closure.
\end{abstract}

\maketitle

\tableofcontents

\section{Introduction}\label{sec:intro}

\subsection{Background}
We study the evolution of two immiscible viscous compressible fluids in the three-dimensional periodic domain $\Omega=\mathbb{T}^3$ for $t\geq 0$. The two phases share a common velocity field and are coupled through an algebraic pressure closure. In Eulerian coordinates, the model reads
\begin{align}\label{eq:TF-orig}
    \begin{dcases}
        \partial_t (\alpha_{\pm} \varrho_{\pm})
        + \div (\alpha_{\pm} \varrho_{\pm} \mathbf{u})
        = 0, \\
         \p_t[ (\alpha_{+} \varrho_{+} +\alpha_{-} \varrho_{-}  )   \mathbf{u} ] 
         + \div [ (\alpha_{+} \varrho_{+} +\alpha_{-} \varrho_{-}  )   \mathbf{u} \otimes \mathbf{u}  ]
         + \nabla p= \mu \Delta \mathbf{u} +(\mu+\lambda) \nabla \div \mathbf{u} 
         , \\
         \alpha_{+} + \alpha_{-}=1, \quad \alpha_{\pm}  \geq 0
         , \\
         p=p_{+}=p_{-}. 
    \end{dcases}
\end{align}
Here $\alpha_\pm=\alpha_\pm(t,x)$ denote the volume fractions of the two phases, $\varrho_\pm=\varrho_\pm(t,x)$ their mass densities, $\mathbf{u} =\mathbf{u} (t,x)$ the common velocity field, and $p$ the common pressure. The constants $\mu$ and $\lambda$ are the shear and bulk viscosity coefficients, respectively.

Following Bresch et al.~\cite{Bre-Muc-Zat-19}, it is convenient to introduce the conservative variables
\begin{align}
 R=\alpha_{+} \varrho_{+},\quad 
 Q=\alpha_{-} \varrho_{-} , \quad 
 Z =  \varrho_{+}.
\end{align}
Then system \eqref{eq:TF-orig} can be rewritten as
\begin{align}\label{eq:TF-R1}
    \begin{dcases}
        \partial_t R + \div (R \mathbf{u} )
        = 0, \\
          \partial_t Q + \div (Q \mathbf{u} )
        = 0, \\
        \p_t [(R+Q) \mathbf{u}] + \div [(R+Q) \mathbf{u} \otimes \mathbf{u}]
        +\nabla p= \mu \Delta \mathbf{u} +(\mu+\lambda) \nabla \div \mathbf{u}.
    \end{dcases}
\end{align}
The algebraic pressure closure means that the phase pressures coincide. Assuming both fluids are isentropic, we have
\begin{align}
p_{+}=( \varrho_{+})^{ \gamma_{+} } =
p_{-}=( \varrho_{-})^{ \gamma_{-} }, \quad \gamma_\pm>1.
\end{align}
As a consequence, the common pressure can be expressed as
\begin{align}
    p(Z)=Z^{\gamma_{+}},
\end{align}
where $Z$ is determined implicitly by $(R,Q)$ through
\begin{align}\label{eq:press-strc}
    \begin{dcases}
       Q=\sp{ 1-\f{R} {Z }  } Z^{\gamma}
       , \quad \gamma=\frac{ \gamma_{+}  }{ \gamma_{-}  }, 
       \\
       R \leq Z. 
    \end{dcases}
\end{align}
Bresch et al.~ \cite{Bre-Muc-Zat-19} showed that \eqref{eq:press-strc} uniquely determines $Z$, and so there exists some function $\mathcal{Z}(\cdot,\cdot)$ such that $Z=\mathcal{Z}(R,Q)$.

The system \eqref{eq:TF-R1}-\eqref{eq:press-strc} has attracted considerable attention in recent years. For semi-stationary flows, Bresch et al.~\cite{Bre-Muc-Zat-19} proved the existence of finite-energy weak solutions in three space dimensions. For viscous compressible two-fluid systems with general pressure laws, including the algebraic closure considered here, Novotn\'y and Pokorn\'y~\cite{Nov-Pok-20} established the existence global finite-energy weak solutions. Li and Zatorska~\cite{Li-Zat-22} obtained a conditional weak--strong uniqueness principle, and this was recently strengthened to an unconditional weak--strong uniqueness result by Li et al.~\cite{Li-Lu-Po-Ew-26} via a suitable relative entropy inequality. In the inviscid case, Li and Zatorska~\cite{Li-Zat-21} proved the existence of infinitely many weak solutions in three dimensions. In one spatial dimension, Li, Sun, and Zatorska~\cite{Li-Sun-Zat-20} studied global weak solutions and their large-time behavior for large initial data. In the class of strong solutions, Piasecki and Zatorska~\cite{Pia-Zat-22}  proved local well-posedness for large initial data and global well-posedness for small data.

It is also useful to place \eqref{eq:TF-R1}-\eqref{eq:press-strc} within the broader literature on compressible two-fluid models. A first important example is the fluid-particle model
\begin{align}\label{eq:TF-R2}
    \begin{dcases}
        \partial_t R + \div (R \mathbf{u} )
        = 0, \\
          \partial_t Q + \div (Q \mathbf{u} )
        = 0, \\
        \p_t [(R+Q) \mathbf{u}] + \div [(R+Q) \mathbf{u} \otimes \mathbf{u}]
        +\nabla P(R,Q)= \mu \Delta \mathbf{u} +(\mu+\lambda) \nabla \div \mathbf{u}, 
    \end{dcases}
\end{align}
with pressure law
\begin{align}\label{pre:R1}
    P(R,Q)= R^{\gamma}+Q^{\beta}, \quad  \gamma,\beta \geq 1. 
\end{align}
For this model, by developing the variable reduction argument, Vasseur, Wen, and Yu~\cite{Vas-Wen-Yu-19} proved the global existence of finite-energy weak solutions under suitable restrictions on $\gamma$ and $\beta$,
namely either 
\begin{align}
 \gamma,\beta > \f{9}{5}, \quad 
 \max\left\{ 
 \f{3\gamma}{4},\gamma-1, \f{3(\gamma+1)}{5}
 \right\}
 <\beta<\min\left\{  
\f{4\gamma}{3}, \gamma+1, \f{5\gamma}{3}-1
\right\},
\end{align}
or
\begin{align}
\gamma>\f{9}{5}, \quad \beta \geq 1, \quad 
\f{1}{c_0} R_0 \leq Q_0 \leq c_0 R_0 \quad 
\text{ for some } c_0 \geq 1. 
\end{align}
Wen~\cite{Wen-21} later improved this result to $\gamma, \beta \geq \f{9}{5}$, with all constraints between adiabatic constants or two densities removed. His result implies, in particular, that  transition to each single-phase flow is allowed.
We also mention the existence theory for dissipative weak solutions due to Li and She~\cite{Li-She-22}, allowing to obtain global-in-time existence result for all finite energy initial data and all $\beta,\gamma>1$. 
Finally, in the framework of finite energy weak solutions,  the inviscid incompressible limit for general initial data in three-dimensional whole space was justified by Kwon and Li \cite{Kwo-Li-19}.

\medskip

A second representative example is the liquid-gas model
\begin{align}\label{eq:TF-R3}
    \begin{dcases}
        \partial_t R + \div (R \mathbf{u} )
        = 0, \\
          \partial_t Q + \div (Q \mathbf{u} )
        = 0, \\
        \p_t (R \mathbf{u}) + \div (R \mathbf{u} \otimes \mathbf{u})
        +\nabla P(R,Q)= \mu \Delta \mathbf{u} +(\mu+\lambda) \nabla \div \mathbf{u}, 
    \end{dcases}
\end{align}
where
\begin{align}\label{pre:R2}
 P(R,Q)=   
 C \sp{
 -b(R,Q)+ \sqrt{   b^2(R,Q) + c(R,Q)}
 },
\end{align}
with 
\begin{align}
b(R,Q)= k_0-R-a_0 Q, \quad 
c(R,Q)= 4 k_0 a_0 Q,
\end{align}
and $k_0,a_0,C$ are positive  constants. 

This model arises naturally in the description of well and pipe liquid-gas two-phase flows. Global weak solutions were established in one dimension by Evje and Karlsen~\cite{Evj-Kar-08}, and in two dimensions under smallness assumptions by Yao, Zhang, and Zhu~\cite{Yao-Zha-Zhu-10}. In the class of strong solutions, Yao, Zhu, and Zi~\cite{Yao-Zhu-Zi-12} justified the incompressible limit toward the classical incompressible Navier--Stokes equations. We also refer to Wen, Yao, and Zhu~\cite{Wen-Yao-Zhu-18} for a broader overview of related results on the liquid-gas models.

\medskip

A third related class is the one-velocity Baer--Nunziato type system for mixtures of non-interacting compressible fluids:
\begin{align}\label{eq:TF-R4}
    \begin{dcases}
        \partial_t \varrho + \div (\varrho \mathbf{u} )
        = 0, \\
         \partial_t z + \div (z \mathbf{u} )
        = 0, \\
        \p_t \alpha+ \mathbf{u} \cdot \Grad \alpha =0 , \\
        \p_t [ (\varrho +z) \mathbf{u} ]  + \div [  (\varrho +z) \mathbf{u} \otimes \mathbf{u} ]
        +\nabla p(f(\alpha) \varrho, g(\alpha)z ) = \mu \Delta \mathbf{u} +(\mu+\lambda) \nabla \div \mathbf{u}, 
    \end{dcases}
\end{align}
where $f,g:[0,1] \mapsto [0,\infty)$ are given functions. 
Under suitable structural assumptions, Novotn\'y~\cite{Nov-20} proved the existence of global finite-energy weak solutions for an initial-boundary value problem for \eqref{eq:TF-R4}. Weak--strong uniqueness was later obtained by Jin and Novotn\'y~\cite{Jin-Nov-19}; further existence results under general boundary conditions and for dissipative turbulent solutions may be found in \cites{Jin-21, Kra-22}. For a non-isentropic version, Kwon, Novotný, and Cheng~\cite{Kwo-Nov-Che-20} proved weak sequential stability, while Kalousek and Ne\v{c}asov\'{a} \cite{Kal-Nec-26} recently established existence of global weak solutions.

\subsection{Motivation and main results}
Despite the substantial progress described above, the low Mach number limit for the two-fluid model with algebraic pressure closure remains poorly understood. The main difficulty comes from the highly nonlinear and implicit structure of the pressure law: unlike models with an explicit pressure $P(R,Q)$, here the pressure is given by $p(Z)=Z^{\gamma_+}$, where $Z={\mathcal{Z}}(R,Q)$ is defined only implicitly by \eqref{eq:press-strc}. This feature makes both the singular limit and the derivation of quantitative estimates significantly more delicate.

The purpose of this paper is to rigorously justify the low Mach number limit for \eqref{eq:TF-R1}-\eqref{eq:press-strc} in the framework of local-in-time strong solutions, and to obtain explicit convergence rates. To this end, we first introduce the nondimensional scaling
\begin{align}
& R(t,x) \longmapsto  R_{\ep}(\ep t,x ),\quad 
Q(t,x) \longmapsto  Q_{\ep}(\ep t,x ),\quad 
\mathbf{u}(t,x) \longmapsto  \ep \mathbf{u}_{\ep}(\ep t,x ), 
\end{align}
together with
\begin{align}
& \mu \longmapsto \ep \mu_{\ep}, \quad 
\lambda \longmapsto \ep \lambda_{\ep}, \quad 
\mu_{\ep}>0, \quad 2 \mu_{\ep} +3 \lambda_{\ep} \geq 0. 
\end{align}
Then $(R_\varepsilon,Q_\varepsilon,u_\varepsilon)$ formally satisfy
\begin{align}\label{eq:two-fluid}
    \begin{dcases}
        \partial_t R_{\ep} +\div (R_{\ep} \vu) = 0, \\
         \partial_t Q_{\ep} +\div (Q_{\ep} \vu) = 0, \\
        \p_t [ (  R_{\ep} + Q_{\ep} ) \vu]+ \div[ (R_{\ep}+Q_{\ep}) \vu\otimes \vu] + \frac{ \nabla p(Z_{\ep})   }{\ep^2}
        = \mu_{\ep} \Delta \vu      +(\mu_{\ep}+\lambda_{\ep}) \nabla \div \vu, 
    \end{dcases}
\end{align}
where $p(Z_\varepsilon)=Z_\varepsilon^{\gamma_+}$ and $Z_\varepsilon$ is linked to $(R_\varepsilon,Q_\varepsilon)$ through
\begin{align}\label{eq:press}
    \begin{dcases}
       Q_{\ep}=\sp{ 1-\f{R_{\ep}}{Z_{\ep}} } Z_{\ep}^{\gamma}
       , \quad \gamma=\frac{ \gamma_{+}  }{ \gamma_{-}  }, 
       \\
       R_{\ep} \leq Z_{\ep}. 
    \end{dcases}
\end{align}
For simplicity, throughout the paper we assume $\mu_\varepsilon=\mu$ and $\lambda_\varepsilon=\lambda$ are constants satisfying
\[
\mu>0,\qquad 2\mu+3\lambda\ge 0.
\]
Accordingly, we work with the following equivalent system
\begin{align}\label{eq:two-fluid-1}
    \begin{dcases}
        \partial_t R_{\ep} +\div (R_{\ep} \vu) = 0, \\
         \partial_t Q_{\ep} +\div (Q_{\ep} \vu) = 0, \\
        \p_t \vu+ \vu \cdot \Grad \vu + \frac{ \nabla p(Z_{\ep})   }{\ep^2   (R_{\ep} + Q_{\ep}) }
        = \f{\mu}{  R_{\ep} + Q_{\ep} } \Delta \vu      
        +\f{  \mu+\lambda  }{  R_{\ep} + Q_{\ep} } \nabla \div \vu . 
    \end{dcases}
\end{align}

If
\begin{align}
R_{\ep} \rightarrow 1 , \quad 
Q_{\ep} \rightarrow 1 , \quad 
\vu \rightarrow \mathbf{u}   \quad \text{ as } 
\varepsilon \rightarrow 0   ,
\end{align} 
then the expected limit system on $\mathbb{T}^3$ is the incompressible Navier--Stokes equations
\begin{align}\label{eq:limit}
    \begin{dcases}
        \partial_t  \mathbf{u} + \mathbf{u} \cdot\nabla  \mathbf{u}
        +\nabla \Pi= \f{\mu}{2} \Delta \mathbf{u} ,
       \\
        \div \mathbf{u}=0. 
    \end{dcases}
\end{align}   

Our first main result proves that, for well-prepared initial data, system \eqref{eq:two-fluid-1} admits a unique classical solution on a time interval independent of $\varepsilon$, together with uniform estimates sufficient to pass to the limit $\varepsilon\to 0$. In particular, $R_\varepsilon$ and $Q_\varepsilon$ converge strongly to the constant state $1$, while $\mathbf{u}_\varepsilon$ converges to the strong solution of \eqref{eq:limit}. Our second main result provides quantitative convergence rates. More precisely, under an additional smallness assumption on the initial energy, we show that
\[
\|R_\varepsilon-1\|_{H^s}^2+\|Q_\varepsilon-1\|_{H^s}^2 \le C\varepsilon^2,
\qquad
\|\mathbf{u}_\varepsilon-\mathbf{u}\|_{L^2}^2 + \int_0^t \|\mathbf{u}_\varepsilon-\mathbf{u}\|_{H^1}^2\,d\tau \le C\varepsilon,
\]
and
\[
\|\div \mathbf{u}_\varepsilon\|_{H^{s-1}} \le C\varepsilon.
\]

These results appear to be the first quantitative low Mach number limit for the compressible two-fluid model with algebraic pressure closure in the strong-solution setting considered here. Compared with the incompressible-limit result for the liquid-gas model obtained in \cite{Yao-Zhu-Zi-12}, our convergence estimates do not require a uniform bound on the quotient of the initial densities. We also note that a recent preprint by Lebot~\cite{Leb-26} studies a related low Mach number problem for a compressible two-phase flow system with algebraic closure and identifies an inhomogeneous incompressible Navier--Stokes system as the limit for well-prepared data. In contrast, our result yields explicit convergence rates and only assumes $\gamma_\pm>1$.

From the technical viewpoint, the proof combines uniform high-order energy estimates with a relative energy argument. The main obstacle is that the pressure is not an explicit function of $(R,Q)$; instead, one has to exploit the structure of the implicitly defined map $(R,Q)\mapsto {\mathcal{Z}}(R,Q)$ and carefully control its derivatives in the singular regime $\varepsilon\to 0$. This makes the commutator estimates and the closure of the energy method substantially more involved than in single-fluid or explicitly closed two-fluid models.

\medskip

The rest of the paper is organized as follows. In Section~\ref{Sec:2} we state the main results. Section~\ref{Sec:3} is devoted to the uniform estimates and the proof of the existence theorem on an $\varepsilon$-independent time interval. In Section~\ref{Sec:4} we establish the relative energy inequality and derive the convergence rates.

\section{Statement of the main results}\label{Sec:2}
For an integer $s\ge 0$, we introduce the energy functionals
{\small{
\begin{align}
 \mathscr{E}_{s} (R,Q,\mathbf{u})
 & :=
 \f{1}{2} \sum_{|\alpha|\leq s} 
 \int_{\Om} \sp{
 \f{1}{\ep^2}  |  \nabla^{\alpha}(R-1)|^2
 +
 \f{1}{\ep^2}  |  \nabla^{\alpha}(Q-1)|^2
 +
 | \nabla^{\alpha} \mathbf{u} |^2
 } dx, \\
  \mathscr{F}_{s} (R,Q,\mathbf{u})
 & :=
  \f{1}{2} \sum_{|\alpha|\leq s} 
 \int_{\Om} \sp{
 \f{  \gamma_{+} \Z^{\gamma_{+}-1} (\p_{R}\Z )(R,Q)      }{ \ep^2 R  }
|  \nabla^{\alpha}(R-1)|^2
  } dx \\
  &\quad 
  + 
  \f{1}{2} \sum_{|\alpha|\leq s} 
 \int_{\Om} \sp{
 \f{  \gamma_{+} \Z^{\gamma_{+}-1} (\p_{Q}\Z )(R,Q)      }{ \ep^2 Q  }
|  \nabla^{\alpha}(Q-1)|^2 
+ (R+Q) | \nabla^{\alpha} \mathbf{u} |^2
  } dx.
\end{align}      }}

When $R$ and $Q$ remain sufficiently close to $1$, the functionals
$\mathscr{E}_{s} (R,Q,\mathbf{u})$ and $\mathscr{F}_{s} (R,Q,\mathbf{u})$ are equivalent.  The functional $\mathscr{F}_{s}(R,Q,\mathbf{u})$ is the natural energy generated by the symmetrized structure of  system \eqref{eq:two-fluid-1},
while $\mathscr{E}_{s} (R,Q,\mathbf{u})$ is more convenient for stating the
uniform bounds and convergence estimates.

We now state the main results of the paper. The first theorem provides a uniform
existence theory to \eqref{eq:two-fluid-1} on a time interval independent of $\varepsilon$, together with the
convergence of strong solutions to the incompressible limit \eqref{eq:limit}.
\begin{thm}\label{Thm1}
Let
\begin{align}
R_{0,\ep}= 1+ \widetilde{R_{0,\ep}}, \quad 
Q_{0,\ep}= 1+ \widetilde{Q_{0,\ep}}, \quad 
\mathbf{u}_{0,\ep}= \mathbf{u}_0 + \widetilde{\mathbf{u}_{0,\ep}}, 
\end{align}
where $\mathbf{u}_0$ satisfies $\div \mathbf{u}_0=0,\mathbf{u}_0\in H^{s+1}(\Om)$ for some integer $s \geq 4$. Assume further that  
\begin{align}
\n{ \widetilde{R_{0,\ep}}  }_{H^s}
\leq  \ep^2 \delta_0, \quad 
\n{ \widetilde{Q_{0,\ep}}  }_{H^s}
\leq  \ep^2 \delta_0, \quad 
\n{ \widetilde{\mathbf{u}_{0,\ep}}  }_{H^{s+1}}
\leq  \ep \delta_0, \quad 
\end{align}
for some sufficiently small constant $\delta_0>0$ independent of $\varepsilon$. 

Then the following statements hold, provided $\varepsilon>0$ is sufficiently small.

\begin{enumerate}
    \item {
    \textbf{Uniform regularity}.  
    There exist positive constants $T_{\ast}$ and $C>0$, independent of $\varepsilon$,
such that \eqref{eq:two-fluid-1} admits a unique classical solution $(R_{\ep},Q_{\ep},\vu)$ with initial data $(R_{0,\ep},Q_{0,\ep},\mathbf{u}_{0,\ep})$ on the time interval $[0,T_{\ast}]$. Moreover, the solution obeys the uniform estimates
   \begin{align}\label{eq:unif-reg}
    \begin{dcases}
        \mathscr{E}_{s} (R_{\ep},Q_{\ep},\vu)(t) 
        +
        \int_0^t \sp{
        \mu \n{ \nabla \vu }_{H^s}^2
        + (\mu+\lambda) \n{ \div \vu }_{H^s}^2
        } d\tau \leq C
        ,
       \\
        \mathscr{E}_{s-1} (\p_t R_{\ep}, \p_tQ_{\ep},\p_t\vu)(t) \\
        \qquad \qquad 
        +
        \int_0^t \sp{
        \mu \n{ \nabla \p_t \vu }_{H^{s-1}}^2
        + (\mu+\lambda) \n{ \div \p_t \vu }_{H^{s-1}}^2
        } d\tau \leq C, 
    \end{dcases}
\end{align}
for all $t\in [0,T_{\ast}]$. 
    }
    \item{
   \textbf{Convergence to the incompressible flow}. 
    As $\ep \rightarrow 0$, we have t
    \begin{align}
    \begin{dcases}
        R_{\ep} \rightarrow 1 \quad \text{ strongly in } 
        L^{\infty}(0,T_{\ast}; H^s(\Om)) \cap \text{Lip} ([0,T_{\ast}]; H^{s-1}(\Om))
         ,
       \\
         Q_{\ep} \rightarrow 1 \quad \text{ strongly in } 
        L^{\infty}(0,T_{\ast}; H^s(\Om)) \cap \text{Lip} ([0,T_{\ast}]; H^{s-1}(\Om))
         ,
       \\
        \vu \rightarrow \mathbf{u} \quad \text{ weakly-$\ast$ in } 
        L^{\infty}(0,T_{\ast}; H^s(\Om)) \cap \text{Lip} ([0,T_{\ast}]; H^{s-1}(\Om))
         ,
       \\
        \vu \rightarrow \mathbf{u} \quad \text{ strongly in } 
        C ( [0,T_{\ast}]; H^{s-\delta'}(\Om) )
    \end{dcases}
\end{align}
for some sufficiently small $\delta'>0$. Here, the limit velocity $ \mathbf{u} \in L^{\infty}(0,T_{\ast}; H^s(\Om)) \cap \text{Lip} ([0,T_{\ast}]; H^{s-1}(\Om))$ is the unique classical solution to the incompressible Navier-Stokes equations \eqref{eq:limit} with initial data $\mathbf{u}_0$. 
    }
\end{enumerate}

\end{thm}

\begin{rem}
The regularity assumption $s\ge 4$ is not expected to be optimal. It is imposed here
to simplify the high-order energy estimates needed for the low Mach number limit.
\end{rem}

Our second result gives quantitative convergence rates.
\begin{thm}\label{Thm2}
Let the assumptions of Theorem \ref{Thm1} be satisfied. In addition, suppose that
\begin{align}
    \int_{\Om} 
    \left[
    \f{1}{\gamma_{+}-1} \sp{ \f{R_{0,\varepsilon}}{ \alpha_{0,\varepsilon} }}^{\gamma_{+}  } 
    \alpha_{0,\varepsilon} 
    +
    \f{1}{\gamma_{-}-1} \sp{ \f{Q_{0,\varepsilon}}{ 1-\alpha_{0,\varepsilon} }}^{\gamma_{-}  } 
    \sp{ 1-\alpha_{0,\varepsilon} } 
    \right]
    dx
    \leq C \ep^3, 
\end{align}
where $\alpha_{0,\varepsilon}:=R_{0,\varepsilon}/Z_{0,\varepsilon}$. 
Then, for all $t \in [0,T_{\ast}]$, the following estimates hold:
\begin{align}\label{eq:rate}
\begin{dcases}
\n{R_{\ep}-1}_{H^s}^2 + \n{Q_{\ep}-1}_{H^s}^2
\leq C \ep^2, 
\\
\n{ \vu- \mathbf{u} }_{L^2}^2 
+ 
\int_0^t \n{ \vu- \mathbf{u} }_{H^1}^2   d\tau 
\leq C \ep, \\
\n{ \div \vu}_{H^{s-1}} 
= 
\n{ \div \sp{ \vu- \mathbf{u}}  }_{H^{s-1}} 
\leq C \ep. 
\end{dcases} 
\end{align}
\end{thm}

\begin{rem}
We give some remarks about the main theorems as follows. 
\begin{enumerate}
    \item{
    Observe that Yao et al. \cite{Yao-Zhu-Zi-12} obtained the convergence rates from the viscous liquid-gas model \eqref{eq:TF-R3}-\eqref{pre:R2} to the classical incompressible Navier-Stokes equations by assuming that the quotient of two initial densities is uniformly bounded, i.e., 
\begin{align}\label{eq:assp}
    \sup_{ \Om} \f{ Q_{0,\ep}  }{ R_{0,\ep}  } \leq 1. 
\end{align}
Compared with Yao et al. \cite{Yao-Zhu-Zi-12}, our convergence estimates \eqref{eq:rate} 
do not require any analogous uniform bound on the ratio of the initial densities. 
    }
    \item{
    After completion of this work, we became aware of the preprint \cite{Leb-26}, which studies
the low Mach number limit for a compressible two-phase flow system with algebraic
closure and derives an inhomogeneous incompressible Navier--Stokes system as the
limit for well-prepared initial data. 
    Our Theorem \ref{Thm2} provides an exact convergence rate, while there is no convergence rate in \cite{Leb-26}. 
    Another difference is that Lebot \cite{Leb-26} studied the low Mach number limit under the assumption that $\gamma_{\pm} \geq 2$, while in our paper, we only require $\gamma_{\pm} > 1$. 
    }
\end{enumerate}

\end{rem}

\section{Proof of Theorem \ref{Thm1}}\label{Sec:3}

To simplify the presentation, throughout this section we assume that
\[
\gamma_+<\gamma_-,
\qquad \text{equivalently, } \gamma=\frac{\gamma_+}{\gamma_-}<1.
\]
By symmetry, the same argument also applies in the case \(\gamma_+>\gamma_-\).
Unless otherwise specified, \(C>0\) denotes a generic constant, independent of
\(\varepsilon\), whose value may change from line to line. The proof combines a
uniform high-order energy method with a contraction argument for a suitable linearized
problem, and is partly inspired by the incompressible limit analysis for the Oldroyd--B
model in \cite{Lei-06} and the viscous liquid-gas model in \cite{Yao-Zhu-Zi-12}.

\subsection{Uniform estimates}\label{subs:unif-est}

Let us define $\mathcal{A}_{T_{\ast},\ep  }(R_{0,\ep}, Q_{0,\ep}, \mathbf{u}_{0,\ep},\delta,C_1,C_2)$ as the subset of $L^{\infty}(0,T_{\ast}; H^s(\Om)) \cap \text{Lip} ([0,T_{\ast}]; H^{s-1}(\Om))$ with $s\geq 4$ such that
\begin{align}
    \begin{dcases}
        | R_{\ep}-1|+|Q_{\ep}-1| < \ep \delta,
       \\
       \mathscr{E}_{s} (R_{\ep},Q_{\ep},\vu)(t) 
        +
        \int_0^t \sp{
        \mu \n{ \nabla \vu }_{H^s}^2
        + (\mu+\lambda) \n{ \div \vu }_{H^s}^2
        } d\tau \leq C_1, \\
        \mathscr{E}_{s-1} (\p_t R_{\ep}, \p_tQ_{\ep},\p_t\vu)(t) 
        +
        \int_0^t \sp{
        \mu \n{ \nabla \p_t \vu }_{H^{s-1}}^2
        + (\mu+\lambda) \n{ \div \p_t \vu }_{H^{s-1}}^2
        } d\tau \leq C_2 
    \end{dcases}
\end{align}
for all $t \in [0,T_{\ast}]$. For $(r_{\ep},q_{\ep},\vv) \in \mathcal{A}_{T_{\ast},\ep  }(R_{0,\ep}, Q_{0,\ep}, \mathbf{u}_{0,\ep},\delta,C_1,C_2)$, we define the map
\begin{align}
 \chi:   (r_{\ep},q_{\ep},\vv) \mapsto 
 (R_{\ep},Q_{\ep},\vu), 
\end{align}
where $(R_{\ep},Q_{\ep},\vu)$ is the unique solution to the following linear problem {\small{
\begin{align}\label{eq:lin}
    \begin{dcases}
        \partial_t R_{\ep} + \vv \cdot \Grad R_{\ep}+ r_{\ep} \div \vu= 0, \\
         \partial_t Q_{\ep} + \vv \cdot \Grad Q_{\ep}+ q_{\ep} \div \vu= 0, \\
        \p_t \vu+ \vv \cdot \Grad \vu + 
        \frac{ \gamma_{+} [ \Z^{\gamma_{+}-1} (\p_{R}\Z ) ](r_{\ep},q_{\ep})    }{\ep^2   (r_{\ep} + q_{\ep}) }
        \Grad R_{\ep}  
        +
        \frac{ \gamma_{+} [ \Z^{\gamma_{+}-1} (\p_{Q}\Z ) ](r_{\ep},q_{\ep})    }{\ep^2   (r_{\ep} + q_{\ep}) }
        \Grad Q_{\ep} 
        \\
        \qquad \qquad \qquad \qquad 
        = \f{\mu}{  r_{\ep} + q_{\ep} } \Delta \vu      
        +\f{  \mu+\lambda  }{  r_{\ep} + q_{\ep} } \nabla \div \vu , 
    \end{dcases}
\end{align}  }}
with initial data
\begin{align}
R_{\ep}(0,\cdot)= 1+ \widetilde{R_{0,\ep}}, \quad 
Q_{\ep}(0,\cdot)= 1+ \widetilde{Q_{0,\ep}}, \quad 
\mathbf{u}_{\ep}(0,\cdot)= \mathbf{u}_0 + \widetilde{\mathbf{u}_{0,\ep}}. 
\end{align}
Our aim is to show that, for properly chosen $T_{\ast},\delta,C_1,C_2$ independent of $\ep$, $\chi$ maps $\mathcal{A}_{T_{\ast},\ep  }(R_{0,\ep}, Q_{0,\ep}, \mathbf{u}_{0,\ep},\delta,C_1,C_2)$ into itself and is a contraction in suitable norm.

For clarity of notations, we denote by $D$ the differential operator of space-time. 
Applying the operator $D^{\alpha}$ to \eqref{eq:lin} gives
\begin{align}\label{eq:D-alpha} 
    \begin{dcases}
       \p_t D^{\alpha} R_{\ep} + \vv \cdot \Grad D^{\alpha} R_{\ep}
       + r_{\ep} \div D^{\alpha} \vu=
       \mathcal{N}_1
        ,
       \\
        \p_t D^{\alpha} Q_{\ep} + \vv \cdot \Grad D^{\alpha} Q_{\ep}
       + q_{\ep} \div D^{\alpha} \vu=
       \mathcal{N}_2
        ,
       \\
       \p_t D^{\alpha} \vu + \vv \cdot \Grad  D^{\alpha} \vu
       + \frac{ \gamma_{+} [ \Z^{\gamma_{+}-1} (\p_{R}\Z ) ](r_{\ep},q_{\ep})    }{\ep^2   (r_{\ep} + q_{\ep}) }
        \Grad   D^{\alpha} R_{\ep}  \\
         \quad 
         +
         \frac{ \gamma_{+} [ \Z^{\gamma_{+}-1} (\p_{Q}\Z ) ](r_{\ep},q_{\ep})    }{\ep^2   (r_{\ep} + q_{\ep}) }
        \Grad D^{\alpha} Q_{\ep} 
        =
        \f{\mu}{  r_{\ep} + q_{\ep} } \Delta D^{\alpha} \vu      
        +\f{  \mu+\lambda  }{  r_{\ep} + q_{\ep} } \nabla \div  D^{\alpha} \vu
        + \mathcal{N}_3, 
    \end{dcases}
\end{align}
where 
\begin{align}
 \mathcal{N}_1
& :=  
- \Big( D^{\alpha} ( \vv \cdot \Grad R_{\ep} ) - \vv \cdot \Grad 
D^{\alpha} R_{\ep}  \Big)
-
\Big( 
D^{\alpha}  ( r_{\ep} \div \vu )
- r_{\ep} \div D^{\alpha} \vu
\Big)
, 
\\
\mathcal{N}_2
& := 
- \Big( D^{\alpha} ( \vv \cdot \Grad Q_{\ep} ) - \vv \cdot \Grad 
D^{\alpha} Q_{\ep}  \Big)
-
\Big( 
D^{\alpha}  ( q_{\ep} \div \vu )
- q_{\ep} \div D^{\alpha} \vu
\Big)
, 
\\
\mathcal{N}_3
& :=  
-\Big(
D^{\alpha} (  \vv \cdot  \Grad \vu ) -\vv \cdot \Grad D^{\alpha} \vu
\Big) \\
& \quad 
-
\left[
D^{\alpha} 
\left( 
\frac{ \gamma_{+} [ \Z^{\gamma_{+}-1} (\p_{R}\Z ) ](r_{\ep},q_{\ep})    }{\ep^2   (r_{\ep} + q_{\ep}) }
        \Grad R_{\ep}  
\right)
- 
\frac{ \gamma_{+} [ \Z^{\gamma_{+}-1} (\p_{R}\Z ) ](r_{\ep},q_{\ep})    }{\ep^2   (r_{\ep} + q_{\ep}) }
        \Grad   D^{\alpha} R_{\ep} 
\right]  \\
\\
& \quad 
-
\left[
D^{\alpha} 
\left( 
\frac{ \gamma_{+} [ \Z^{\gamma_{+}-1} (\p_{Q}\Z ) ](r_{\ep},q_{\ep})    }{\ep^2   (r_{\ep} + q_{\ep}) }
        \Grad Q_{\ep}  
\right)
- 
\frac{ \gamma_{+} [ \Z^{\gamma_{+}-1} (\p_{Q}\Z ) ](r_{\ep},q_{\ep})    }{\ep^2   (r_{\ep} + q_{\ep}) }
        \Grad   D^{\alpha} Q_{\ep} 
\right] \\
& \quad 
+
\mu 
\left[ 
D^{\alpha} \left(
\f{\Delta \vu}{ r_{\ep}+q_{\ep} }
\right)- 
\f{1}{r_{\ep}+q_{\ep}} \Delta D^{\alpha} \vu
\right] \\
&
\quad 
+
(\mu +\lambda)
\left[ 
D^{\alpha} \left(
\f{\Grad \div \vu}{ r_{\ep}+q_{\ep} }
\right)- 
\f{1}{r_{\ep}+q_{\ep}} \Grad \div D^{\alpha} \vu
\right] .
\end{align}
To proceed, we take $L^2$-inner product of \eqref{eq:D-alpha}$_1$ with $\frac{ \gamma_{+} [ \Z^{\gamma_{+}-1} (\p_{R}\Z ) ](r_{\ep},q_{\ep})    }{\ep^2   r_{\ep}  } D^{\alpha}(R_{\ep}-1)$, \eqref{eq:D-alpha}$_2$ with $\frac{ \gamma_{+} [ \Z^{\gamma_{+}-1} (\p_{Q}\Z ) ](r_{\ep},q_{\ep})    }{\ep^2    q_{\ep} } D^{\alpha}(Q_{\ep}-1)$ and\eqref{eq:D-alpha}$_3$ with $(r_{\ep}+q_{\ep}) D^{\alpha} \vu$. Then, we add them together, making use of integration by parts and some cancellations, to arrive at 
\begin{align}
& \f{1}{2} \f{d}{dt} 
\int_{\Om} 
\Big[
\frac{ \gamma_{+} [ \Z^{\gamma_{+}-1} (\p_{R}\Z ) ](r_{\ep},q_{\ep})    }{\ep^2   r_{\ep}  } | D^{\alpha}(R_{\ep}-1)|^2 \\
&  \qquad \qquad +
\frac{ \gamma_{+} [ \Z^{\gamma_{+}-1} (\p_{Q}\Z ) ](r_{\ep},q_{\ep})    }{\ep^2    q_{\ep} } |D^{\alpha}(Q_{\ep}-1)|^2
+ 
(r_{\ep}+q_{\ep}) |D^{\alpha} \vu|^2
\Big] dx \\
& \qquad \qquad +
\mu \n{ \Grad D^{\alpha}\vu  }_{L^2}^2 
+ (\mu+\lambda) \n{ \div D^{\alpha}\vu  }_{L^2}^2 \\
& \qquad 
= 
\sum_{\ell=1}^{7} I_{\ell}, 
\end{align}
where 
\begin{align}
I_{1} & :=     
\f{1}{2}  
\int_{\Om} 
\Big[
\f{1}{\ep^2} | D^{\alpha}(R_{\ep}-1)|^2
\p_t \sp{ 
\frac{ \gamma_{+} [ \Z^{\gamma_{+}-1} (\p_{R}\Z ) ](r_{\ep},q_{\ep})    }{  r_{\ep}  }
}     \\
& \qquad + 
\f{1}{\ep^2} | D^{\alpha}(Q_{\ep}-1)|^2
\p_t \sp{ 
\frac{ \gamma_{+} [ \Z^{\gamma_{+}-1} (\p_{Q}\Z ) ](r_{\ep},q_{\ep})    }{  q_{\ep}  }
} 
+
 |D^{\alpha} \vu|^2   \p_t (r_{\ep}+q_{\ep}) \Big] dx
, \\ 
I_{2} & :=    
\f{1}{2}  
\int_{\Om}  \Big[
| D^{\alpha} (R_{\ep}-1) |^2 
\div \sp{ \frac{ \gamma_{+} [ \Z^{\gamma_{+}-1} (\p_{R}\Z ) ](r_{\ep},q_{\ep})    }{\ep^2   r_{\ep}  }  \vv   }  \\
& \qquad \qquad 
+ 
| D^{\alpha} (Q_{\ep}-1) |^2 
\div \sp{ \frac{ \gamma_{+} [ \Z^{\gamma_{+}-1} (\p_{Q}\Z ) ](r_{\ep},q_{\ep})    }{\ep^2   q_{\ep}  }  \vv   }  \\
& \qquad \qquad 
+
| D^{\alpha} \vu |^2  \div \sp{ (r_{\ep}+q_{\ep}) \vv   } 
\Big] dx
, \\ 
I_{3} & :=   \int_{\Om}   
\Big[ 
D^{\alpha} (R_{\ep}-1)  D^{\alpha}
\vu \cdot \Grad \frac{ \gamma_{+} [ \Z^{\gamma_{+}-1} (\p_{R}\Z ) ](r_{\ep},q_{\ep})    }{\ep^2     } \\
& \qquad \qquad 
+
D^{\alpha} (Q_{\ep}-1)  D^{\alpha}
\vu \cdot \Grad \frac{ \gamma_{+} [ \Z^{\gamma_{+}-1} (\p_{Q}\Z ) ](r_{\ep},q_{\ep})    }{\ep^2     }  \Big] dx
, \\ 
I_{4} & :=   - \int_{\Om}   
\left[
\Big( D^{\alpha} ( \vv \cdot \Grad R_{\ep} ) - \vv \cdot \Grad 
D^{\alpha} R_{\ep}  \Big)
+
\Big( 
D^{\alpha}  ( r_{\ep} \div \vu )
- r_{\ep} \div D^{\alpha} \vu
\Big) \right] \\
& \qquad  \qquad \qquad  \qquad \qquad 
\times \frac{ \gamma_{+} [ \Z^{\gamma_{+}-1} (\p_{R}\Z ) ](r_{\ep},q_{\ep})    }{\ep^2   r_{\ep}  } D^{\alpha}(R_{\ep}-1)  
dx  \\
& \qquad 
- \int_{\Om} \left[
\Big( D^{\alpha} ( \vv \cdot \Grad Q_{\ep} ) - \vv \cdot \Grad 
D^{\alpha} Q_{\ep}  \Big)
+
\Big( 
D^{\alpha}  ( q_{\ep} \div \vu )
- q_{\ep} \div D^{\alpha} \vu
\Big) \right] \\
& \qquad  \qquad \qquad  \qquad \qquad 
\times \frac{ \gamma_{+} [ \Z^{\gamma_{+}-1} (\p_{Q}\Z ) ](r_{\ep},q_{\ep})    }{\ep^2    q_{\ep} } D^{\alpha}(Q_{\ep}-1)
dx
, \\ 
I_{5} & :=   - \int_{\Om}   
(r_{\ep}+q_{\ep}) 
\Big(
D^{\alpha} (  \vv \cdot  \Grad \vu ) -\vv \cdot \Grad D^{\alpha} \vu
\Big) \cdot D^{\alpha}  \vu dx
, \\ 
I_{6} & :=  - \int_{\Om}  (r_{\ep}+q_{\ep})  
\Big[
D^{\alpha} 
\left( 
\frac{ \gamma_{+} [ \Z^{\gamma_{+}-1} (\p_{R}\Z ) ](r_{\ep},q_{\ep})    }{\ep^2   (r_{\ep} + q_{\ep}) }
        \Grad R_{\ep}  
\right)   \\
& \qquad  \qquad \qquad  \qquad \qquad  
- 
\frac{ \gamma_{+} [ \Z^{\gamma_{+}-1} (\p_{R}\Z ) ](r_{\ep},q_{\ep})    }{\ep^2   (r_{\ep} + q_{\ep}) }
        \Grad   D^{\alpha} R_{\ep} 
\Big]  \cdot D^{\alpha}  \vu dx  \\
& \qquad 
- \int_{\Om}  (r_{\ep}+q_{\ep})
\Big[
D^{\alpha} 
\left( 
\frac{ \gamma_{+} [ \Z^{\gamma_{+}-1} (\p_{Q}\Z ) ](r_{\ep},q_{\ep})    }{\ep^2   (r_{\ep} + q_{\ep}) }
        \Grad Q_{\ep}  
\right)  \\
& \qquad  \qquad \qquad  \qquad \qquad 
- 
\frac{ \gamma_{+} [ \Z^{\gamma_{+}-1} (\p_{Q}\Z ) ](r_{\ep},q_{\ep})    }{\ep^2   (r_{\ep} + q_{\ep}) }
        \Grad   D^{\alpha} Q_{\ep} 
\Big] \cdot D^{\alpha}  \vu dx
,       \\ 
I_{7} & :=   
\mu \int_{\Om}  (r_{\ep}+q_{\ep})
\left[ 
D^{\alpha} \left(
\f{\Delta \vu}{ r_{\ep}+q_{\ep} }
\right)- 
\f{1}{r_{\ep}+q_{\ep}} \Delta D^{\alpha} \vu
\right] \cdot D^{\alpha}  \vu dx \\
& \qquad 
+ (\mu+\lambda) \int_{\Om}  (r_{\ep}+q_{\ep})
\left[ 
D^{\alpha} \left(
\f{\Grad \div \vu}{ r_{\ep}+q_{\ep} }
\right)- 
\f{1}{r_{\ep}+q_{\ep}} \Grad \div D^{\alpha} \vu
\right] \cdot D^{\alpha}  \vu dx
.  
\end{align}
The key issue is to estimate $I_{\ell},\ell=1,...,7$ suitably.

Before turning to the detailed estimates, we give some basic computations as follows.  
\begin{align}\label{deriv-pre}
\p_{R} p(\Z)  &= \gamma_{+} \Z^{\gamma_{+} -1 } \p_{R} \Z, \quad 
\p_{Q} p(\Z) = \gamma_{+} \Z^{\gamma_{+} -1 } \p_{Q} \Z, \\
\p_{RR} p(\Z)   & =
\gamma_{+}  (\gamma_{+} -1) \Z^{\gamma_{+} -2 } (\p_{R} \Z)^2
+ \gamma_{+} \Z^{\gamma_{+} -1 } \p_{RR} \Z, \\
\p_{RQ} p(\Z) & = \gamma_{+}  (\gamma_{+} -1) \Z^{\gamma_{+} -2 } \p_{Q} \Z \p_{R} \Z + \gamma_{+} \Z^{\gamma_{+} -1 } \p_{RQ} \Z, \\
\p_{R} \Z &=
\f{ \Z^{\gamma-1} }
{
\gamma \Z^{\gamma-1} - R (\gamma-1) \Z^{\gamma-2}
} , \quad 
\p_{Q} \Z=
\f{ 1 }
{
\gamma \Z^{\gamma-1} - R (\gamma-1) \Z^{\gamma-2}
} , \\
\p_{RR} \Z &=
\f{ 
\left[(\gamma-1) \Z^{\gamma-2} \p_{R} \Z \right]
\left[ \gamma \Z^{\gamma-1} - R (\gamma-1) \Z^{\gamma-2} \right]
}
{
\left[ \gamma \Z^{\gamma-1} - R (\gamma-1) \Z^{\gamma-2} \right]^2
} \\
& \quad 
- 
\f{
\Z^{\gamma-1} 
\left[  
\gamma (\gamma-1) \Z^{\gamma-2} \p_{R} \Z 
- (\gamma-1) \Z^{\gamma-2} - (\gamma-1) R (\gamma-2) \Z^{\gamma-3}  \p_{R} \Z 
\right]
}
{
\left[ \gamma \Z^{\gamma-1} - R (\gamma-1) \Z^{\gamma-2} \right]^2
}, \\
\p_{RQ} \Z &=
\f{
\left[(\gamma-1) \Z^{\gamma-2} \p_{Q} \Z \right]
\left[ \gamma \Z^{\gamma-1} - R (\gamma-1) \Z^{\gamma-2} \right]
}
{
\left[ \gamma \Z^{\gamma-1} - R (\gamma-1) \Z^{\gamma-2} \right]^2
} \\
& \quad 
- 
\f{
\Z^{\gamma-1} 
\left[  
\gamma (\gamma-1) \Z^{\gamma-2} \p_{Q} \Z 
 - R (\gamma-1)  (\gamma-2) \Z^{\gamma-3}  \p_{Q} \Z 
\right]
}
{
\left[ \gamma \Z^{\gamma-1} - R (\gamma-1) \Z^{\gamma-2} \right]^2
}, \\
\p_{QQ} \Z &= - 
\f{
\gamma (\gamma-1) \Z^{\gamma-2} \p_{Q} \Z 
- R (\gamma-1) (\gamma-2) \Z^{\gamma-3}  \p_{Q} \Z
}
{
\left[ \gamma \Z^{\gamma-1} - R (\gamma-1) \Z^{\gamma-2} \right]^2
}.
\end{align}
Next, we recall an observation from \cites{Li-Sun-Zat-20,Li-Zat-21}. Suppose that 
\begin{align}
0<\underline{R} \leq R(t,x) \leq \overline{R} <\infty, \quad
0<\underline{Q} \leq Q(t,x) \leq \overline{Q} <\infty,
\end{align}
for all $[0,T_{\ast}] \times \Om$. Then the structure conditions \eqref{eq:press} ensure that there exist positive constants $\underline{Z},\overline{Z}$ such that 
\begin{align}\label{bound-Z}
   0<\underline{Z} \leq \Z(t,x) \leq \overline{Z} <\infty, 
\end{align}
for all $[0,T_{\ast}] \times \Om$. Indeed, as it was shown in \cite{Li-Zat-21}, it suffices to take
\begin{align}
  \underline{Z}=  \underline{R} ,\quad 
  \overline{Z}=
  \max\left\{ 2\overline{R}, \sp{2\overline{Q}}^{1/\gamma} \right\}.
\end{align}


\noindent 
{\it Step 1. The estimates for $I_1-I_3$.} 
To begin with, we choose $\delta>0$ sufficiently small such that $|r_{\ep}-1|\leq 1/2$ and $|q_{\ep}-1|\leq 1/2$ hold, which means $1/2 \leq r_{\ep},q_{\ep} \leq 3/2$. It follows from the Sobolev's embedding inequality that 
\begin{align}
|I_1| & \leq C \n{ \f{ D^{\alpha}(R_{\ep}-1)  }{\ep}  }_{L^2}^2
   \n{ 
   \f{
   \left[
   \p_{RR}p(\Z)     
    \p_t r_{\ep} + \p_{RQ} p(\Z)     
    \p_t q_{\ep} \right] r_{\ep} -  \p_{R}p(\Z) \p_t r_{\ep}
   }
   {
   r_{\ep}^2
   }  (r_{\ep},q_{\ep}) 
   }_{L^{\infty}} \\
   & \quad 
   + 
   C \n{ \f{ D^{\alpha}(Q_{\ep}-1)  }{\ep}  }_{L^2}^2
   \n{ 
   \f{
   \left[
   \p_{QR}p(\Z)     
    \p_t r_{\ep} + \p_{QQ} p(\Z)     
    \p_t q_{\ep} \right] q_{\ep} -  \p_{Q}p(\Z) \p_t q_{\ep}
   }
   {
   q_{\ep}^2
   }  (r_{\ep},q_{\ep}) 
   }_{L^{\infty}} \\
   & \quad 
   + 
   \n{D^{\alpha} \vu }_{L^2}^2 
   \n{  \p_t r_{\ep} +\p_t q_{\ep}  }_{L^{\infty}}\\
   & 
   \leq 
   C \sp{
   \n{\p_t  r_{\ep}}_{H^{s-2}} + \n{\p_t  q_{\ep}}_{H^{s-2}}
   }
   \sp{
   \n{ \f{ D^{\alpha}(R_{\ep}-1)  }{\ep}  }_{L^2}^2
   +\n{ \f{ D^{\alpha}(Q_{\ep}-1)  }{\ep}  }_{L^2}^2
   +\n{D^{\alpha} \vu }_{L^2}^2 
   }\\
   & 
   \leq 
   C  C_2^{\f{1}{2}} \ep 
   \sp{
   \n{ \f{ D^{\alpha}(R_{\ep}-1)  }{\ep}  }_{L^2}^2
   +\n{ \f{ D^{\alpha}(Q_{\ep}-1)  }{\ep}  }_{L^2}^2
   +\n{D^{\alpha} \vu }_{L^2}^2 
   }.
\end{align}
Here, we have used the fact that\footnote{This fact will be extensively used in the sequel, without pointing it out each time.} $(r_{\ep},q_{\ep},\vv) \in \mathcal{A}_{T_{\ast},\ep  }(R_{0,\ep}, Q_{0,\ep}, \mathbf{u}_{0,\ep},\delta,C_1,C_2)$ and the uniform bounds
\begin{align}\label{unif-bd}
 & \n{ \p_{R}p(\Z)(r_{\ep},q_{\ep})  }_{L^{\infty}} \leq C, 
\quad 
\n{ \p_{Q}p(\Z)(r_{\ep},q_{\ep})  }_{L^{\infty}} \leq C, \\
&
\n{ \p_{RR}p(\Z)(r_{\ep},q_{\ep})  }_{L^{\infty}} \leq C, 
\quad 
\n{ \p_{RQ}p(\Z)(r_{\ep},q_{\ep})  }_{L^{\infty}} \leq C, \\
& \n{ \p_{QQ}p(\Z)(r_{\ep},q_{\ep})  }_{L^{\infty}} \leq C, 
\quad 0<C^{-1} \leq r_{\ep}, q_{\ep}, r_{\ep}^{-1}, q_{\ep}^{-1} \leq C,
\end{align}
which follows from \eqref{deriv-pre}, \eqref{bound-Z} and the two-sided bounds of $r_{\ep}, q_{\ep}$. In the same spirit, we see that 
\begin{align}
  |I_2| & \leq C \n{ \f{ D^{\alpha}(R_{\ep}-1)  }{\ep}  }_{L^2}^2  
  \Bigg[  
  \n{\nabla \vv}_{L^{\infty}} 
  \n{   \frac{ \gamma_{+} [ \Z^{\gamma_{+}-1} (\p_{R}\Z ) ](r_{\ep},q_{\ep})    }{  r_{\ep}  }   }_{L^{\infty}} \\
  & \qquad 
  + \n{ \vv}_{L^{\infty}} 
  \n{
  \f{
  \left[
   \p_{RR}p(\Z)     
    \nabla r_{\ep} + \p_{RQ} p(\Z)     
    \nabla q_{\ep} \right] r_{\ep} -  \p_{R}p(\Z) \nabla r_{\ep}
  }{r_{\ep}^2} (r_{\ep},q_{\ep}) 
   }_{L^{\infty}}  \Bigg] \\
& \quad 
+ 
C \n{ \f{ D^{\alpha}(Q_{\ep}-1)  }{\ep}  }_{L^2}^2   
  \Bigg[  
  \n{\nabla \vv}_{L^{\infty}} 
  \n{   \frac{ \gamma_{+} [ \Z^{\gamma_{+}-1} (\p_{Q}\Z ) ](r_{\ep},q_{\ep})    }{  q_{\ep}  }   }_{L^{\infty}} \\
  & \qquad 
  + \n{ \vv}_{L^{\infty}} 
  \n{
  \f{
  \left[
   \p_{RQ}p(\Z)     
    \nabla r_{\ep} + \p_{QQ} p(\Z)     
    \nabla q_{\ep} \right] q_{\ep} -  \p_{Q}p(\Z) \nabla q_{\ep}
  }{q_{\ep}^2} (r_{\ep},q_{\ep}) 
   }_{L^{\infty}} \Bigg] \\
& \quad 
+ C \n{ D^{\alpha} \vu }_{L^2}^2
\Big[ 
\n{\vv}_{L^{\infty}}  \sp{ \n{\nabla r_{\ep}}_{L^{\infty}}+\n{\nabla q_{\ep}}_{L^{\infty}}   }  \\
& \qquad \qquad \qquad \qquad
+ 
\sp{ \n{ r_{\ep}}_{L^{\infty}}+\n{ q_{\ep}}_{L^{\infty}}   }
\n{ \div \vv}_{L^{\infty}} 
\Big] \\
& \leq C \Big[
\n{\vv}_{L^{\infty}}  \sp{ \n{\nabla r_{\ep}}_{L^{\infty}}+\n{\nabla q_{\ep}}_{L^{\infty}}   } 
+ \n{\nabla \vv}_{L^{\infty}}
\Big] \\
& \qquad \times 
\left[
\n{ \f{ D^{\alpha}(R_{\ep}-1)  }{\ep}  }_{L^2}^2 
+\n{ \f{ D^{\alpha}(Q_{\ep}-1)  }{\ep}  }_{L^2}^2 
+\n{ D^{\alpha} \vu }_{L^2}^2
\right] \\
& \leq C
\left[ 
\n{\vv}_{H^{s-2}} \ep \sp{ 
\n{ \f{\nabla r_{\ep}}{\ep}  }_{H^{s-2}}   
+ \n{ \f{\nabla q_{\ep}}{\ep}  }_{H^{s-2}} 
}
+ \n{\nabla \vv}_{H^{s-2}}
\right]
\\
& \qquad \times 
\left[
\n{ \f{ D^{\alpha}(R_{\ep}-1)  }{\ep}  }_{L^2}^2 
+\n{ \f{ D^{\alpha}(Q_{\ep}-1)  }{\ep}  }_{L^2}^2 
+\n{ D^{\alpha} \vu }_{L^2}^2
\right] \\
& \leq C 
\sp{\ep C_1 + C_1^{\f{1}{2}}   } 
\left[
\n{ \f{ D^{\alpha}(R_{\ep}-1)  }{\ep}  }_{L^2}^2 
+\n{ \f{ D^{\alpha}(Q_{\ep}-1)  }{\ep}  }_{L^2}^2 
+\n{ D^{\alpha} \vu }_{L^2}^2
\right].
\end{align}
Next, we estimate $I_3$ as 
\begin{align}
  |I_3| & \leq \n{ \f{ D^{\alpha}(R_{\ep}-1)  }{\ep}  }_{L^2}
  \n{ D^{\alpha} \vu }_{L^2} 
  \n{ 
  \f{
  \p_{RR}p(\Z)     
    \nabla r_{\ep} + \p_{RQ} p(\Z)     
    \nabla q_{\ep}   
  }{\ep} (r_{\ep},q_{\ep})
  }_{L^{\infty}} \\
  & \quad 
  + 
  \n{ \f{ D^{\alpha}(Q_{\ep}-1)  }{\ep}  }_{L^2}
  \n{ D^{\alpha} \vu }_{L^2} 
  \n{ 
  \f{
  \p_{QR}p(\Z)     
    \nabla r_{\ep} + \p_{QQ} p(\Z)     
    \nabla q_{\ep}   
  }{\ep} (r_{\ep},q_{\ep})
  }_{L^{\infty}} \\
  & \leq 
  C 
  \sp{ 
  \n{\f{\nabla r_{\ep}}{\ep}}_{L^{\infty}}
  +
  \n{\f{\nabla q_{\ep}}{\ep}}_{L^{\infty}}
  } 
  \left[
\n{ \f{ D^{\alpha}(R_{\ep}-1)  }{\ep}  }_{L^2}^2 
+\n{ \f{ D^{\alpha}(Q_{\ep}-1)  }{\ep}  }_{L^2}^2 
+\n{ D^{\alpha} \vu }_{L^2}^2
\right] \\
& \leq 
C \sp{ 
\n{ \f{\nabla r_{\ep}}{\ep}  }_{H^{s-2}}   
+ \n{ \f{\nabla q_{\ep}}{\ep}  }_{H^{s-2}} 
}  \left[
\n{ \f{ D^{\alpha}(R_{\ep}-1)  }{\ep}  }_{L^2}^2 
+\n{ \f{ D^{\alpha}(Q_{\ep}-1)  }{\ep}  }_{L^2}^2 
+\n{ D^{\alpha} \vu }_{L^2}^2
\right] \\
& \leq C C_1^{\f{1}{2}}
\left[
\n{ \f{ D^{\alpha}(R_{\ep}-1)  }{\ep}  }_{L^2}^2 
+\n{ \f{ D^{\alpha}(Q_{\ep}-1)  }{\ep}  }_{L^2}^2 
+\n{ D^{\alpha} \vu }_{L^2}^2
\right] . 
\end{align}
 
\noindent 
{\it Step 2. The estimates for $I_4-I_7$ when $D^{\alpha}=\p_t$.} 
By the uniform bounds \eqref{unif-bd} and the Sobolev's embedding inequality, we see that 
\begin{align}
  |I_4| & \leq   
  C \n{ \ep^{-1}\p_t R_{\ep} }_{L^2}
   \sp{
   \n{\p_t \vv}_{L^{\infty}} \n{ \ep^{-1}\nabla R_{\ep} }_{L^2}
   + \n{ \div \vu}_{L^2} \n{ \ep^{-1} \p_t r_{\ep} }_{L^{\infty}}
   } \\
   & \quad 
   + 
   C \n{ \ep^{-1}\p_t Q_{\ep} }_{L^2}
   \sp{
   \n{\p_t \vv}_{L^{\infty}} \n{ \ep^{-1}\nabla Q_{\ep} }_{L^2}
   + \n{ \div \vu}_{L^2} \n{ \ep^{-1} \p_t q_{\ep} }_{L^{\infty}}
   } \\
   & \leq 
   C \sp{
   \n{ \ep^{-1}\nabla R_{\ep} }_{L^2}^2 + \n{ \ep^{-1}\nabla Q_{\ep} }_{L^2}^2+ \n{\nabla \vu}_{L^2}^2
   } \\
   & \quad 
   + C 
   \sp{
   \n{ \ep^{-1}\p_t R_{\ep} }_{L^2}^2
   +
   \n{ \ep^{-1}\p_t Q_{\ep} }_{L^2}^2
   }
   \sp{
   \n{\p_t \vv}_{L^{\infty}}^2 
   + \n{ \ep^{-1} \p_t r_{\ep} }_{L^{\infty}}^2
   +\n{ \ep^{-1} \p_t q_{\ep} }_{L^{\infty}}^2
   } \\
   & \leq 
   C \sp{
   \n{ \ep^{-1}\nabla R_{\ep} }_{L^2}^2 + \n{ \ep^{-1}\nabla Q_{\ep} }_{L^2}^2+ \n{\nabla \vu}_{L^2}^2
   } \\
   & \quad 
   +
   C C_2
   \sp{
   \n{ \ep^{-1}\p_t R_{\ep} }_{L^2}^2
   +
   \n{ \ep^{-1}\p_t Q_{\ep} }_{L^2}^2
   } , \\
  |I_5| & \leq   
  C \n{ \p_t \vu}_{L^2} \n{ \nabla \vu}_{L^2}  
  \n{\p_t \vv}_{L^{\infty}} 
  \sp{
  \n{r_{\ep}}_{L^{\infty}}+\n{q_{\ep}}_{L^{\infty}} 
  } \\
  & 
  \leq C \n{ \nabla \vu}_{L^2}^2
  + C C_2 \n{ \p_t \vu}_{L^2}^2, \\
   |I_6| & \leq   
   \f{C}{\ep^2} 
   \n{
   \p_t \sp{
   \f{  [\p_{R}p(\Z)] (r_{\ep},q_{\ep}) }{r_{\ep}+q_{\ep}}
   } \nabla R_{\ep}
   }_{L^2} \n{ \p_t \vu}_{L^2} \\
   & \quad 
   + 
   \f{C}{\ep^2} 
   \n{
   \p_t \sp{
   \f{  [\p_{Q}p(\Z)] (r_{\ep},q_{\ep}) }{r_{\ep}+q_{\ep}}
   } \nabla Q_{\ep}
   }_{L^2} \n{ \p_t \vu}_{L^2} \\
   & \leq 
   C \sp{
   \n{\ep^{-1} \p_t r_{\ep} }_{L^{\infty}}
   + \n{\ep^{-1} \p_t q_{\ep} }_{L^{\infty}}
   } 
   \sp{
   \n{ \ep^{-1}\nabla R_{\ep} }_{L^2}
   +
   \n{ \ep^{-1}\nabla Q_{\ep} }_{L^2}
   }   \n{ \p_t \vu}_{L^2} \\
   & \leq 
   C \sp{
   \n{ \ep^{-1}\nabla R_{\ep} }_{L^2}^2
   +
   \n{ \ep^{-1}\nabla Q_{\ep} }_{L^2}^2
   } 
   + C C_2 \n{ \p_t \vu}_{L^2}^2, \\
   |I_7| & \leq   
   C\left(
   \n{
   \p_t \sp{ \f{1}{r_{\ep}+q_{\ep}} } 
   \Delta \vu 
   } 
   \n{ \p_t \vu}_{L^2}
   +
   \n{
   \p_t \sp{ \f{1}{r_{\ep}+q_{\ep}} } 
   \nabla \div \vu 
   } 
   \n{ \p_t \vu}_{L^2}
   \right) \\
   & \leq 
   C \ep
   \sp{
   \n{\ep^{-1} \p_t r_{\ep} }_{L^{\infty}}
   + \n{\ep^{-1} \p_t q_{\ep} }_{L^{\infty}}
   } 
   \sp{
   \n{\Delta \vu}_{L^2} + \n{\nabla \div \vu}_{L^2}
   }\n{ \p_t \vu}_{L^2}
   \\
   & \leq 
   \f{\mu}{100} \n{\Delta \vu}_{L^2}^2 + 
   \f{\mu+\lambda}{100} \n{\nabla \div \vu}_{L^2}^2 
   + C \ep^2  C_2  \n{ \p_t \vu}_{L^2}^2. 
\end{align}

\noindent 
{\it Step 3. The estimates for $I_4-I_7$ when $D^{\alpha}=\nabla^{\alpha},|\alpha|\leq s$.} 
It follows from \eqref{unif-bd}, Sobolev's embedding inequality and Lemma \ref{lemm:prod} that 
\begin{align}
|I_4| & \leq C 
\Big( 
\n{\nabla \vv}_{L^{\infty}} \n{\nabla^s R_{\ep}}_{L^2}
+ \n{\nabla R_{\ep}}_{L^{\infty}} \n{\nabla^s \vv}_{L^2} \\
& \quad 
+ \n{\nabla r_{\ep}}_{L^{\infty}} \n{\nabla^{s-1} \div \vu }_{L^2}
+ \n{  \div \vu}_{L^{\infty}} + \n{\nabla^s r_{\ep}}_{L^2}
\Big) 
\ep^{-2} \n{ \nabla^{\alpha} (R_{\ep}-1) }_{L^2} \\
& \quad 
+
C 
\Big( 
\n{\nabla \vv}_{L^{\infty}} \n{\nabla^s Q_{\ep}}_{L^2}
+ \n{\nabla Q_{\ep}}_{L^{\infty}} \n{\nabla^s \vv}_{L^2} \\
& \quad 
+ \n{\nabla q_{\ep}}_{L^{\infty}} \n{\nabla^{s-1} \div \vu }_{L^2}
+ \n{  \div \vu}_{L^{\infty}} + \n{\nabla^s q_{\ep}}_{L^2}
\Big) 
\ep^{-2} \n{ \nabla^{\alpha} (Q_{\ep}-1) }_{L^2} \\
& \leq 
C C_1^{\f{1}{2}}
\sp{
\n{ \ep^{-1} (R_{\ep}-1)   }_{H^s}^2 
+ \n{ \ep^{-1} (Q_{\ep}-1)   }_{H^s}^2  
+\n{  \vu }_{H^s}^2
}. 
\end{align}
Similarly, we see that 
\begin{align}
|I_5| & \leq 
C \n{\nabla^{\alpha} \vu}_{L^2}
\sp{
\n{\nabla \vv}_{L^{\infty}} \n{\nabla \vu}_{H^{s-1}}
+ \n{\nabla \vu}_{L^{\infty}} \n{\vv}_{H^s}
} \\
& \leq 
C C_1^{\f{1}{2}}   \n{  \vu }_{H^s}^2. 
\end{align}
By \eqref{unif-bd}, Sobolev's embedding inequality, Lemma \ref{lemm:prod} and Lemma \ref{lemm:compo}, we infer that {\small{
\begin{align}
|I_6| & \leq 
C \Big(
\n{
  \f{
  \left[
   \p_{RR}p(\Z)     
    \nabla r_{\ep} + \p_{RQ} p(\Z)     
    \nabla q_{\ep} \right] (r_{\ep} +q_{\ep})  -  \p_{R}p(\Z) \nabla (r_{\ep}+ q_{\ep})
  }{  (r_{\ep} +q_{\ep}) ^2}     
   }_{L^{\infty}}
   \n{ \nabla^{s-1} \nabla R_{\ep}}_{L^2} \\
   & \quad 
   + \n{ \nabla R_{\ep} }_{L^{\infty}}
   \n{ 
   \nabla^{s} \sp{ 
   \f{\p_{R}p(\Z)}{ r_{\ep} +q_{\ep}  }  
   }   
   }_{L^2}
   \Big) \ep^{-2} \n{  \vu }_{H^s} \\
& \quad 
+ 
C \Big(
\n{
  \f{
  \left[
   \p_{QR}p(\Z)     
    \nabla r_{\ep} + \p_{QQ} p(\Z)     
    \nabla q_{\ep} \right] (r_{\ep} +q_{\ep})  - 
    \p_{Q}p(\Z) \nabla (r_{\ep} + q_{\ep})   
  }{(r_{\ep} +q_{\ep})^2 } 
   }_{L^{\infty}}
   \n{ \nabla^{s-1} \nabla Q_{\ep}}_{L^2} \\
   & \quad 
   + \n{ \nabla Q_{\ep} }_{L^{\infty}}
   \n{ 
   \nabla^{s} \sp{ 
   \f{\p_{Q}p(\Z)}{ r_{\ep} +q_{\ep}) } 
   } 
   }_{L^2}
   \Big) \ep^{-2} \n{  \vu }_{H^s} \\
   & \leq 
   C \n{\vu}_{H^s}
   \sp{
  \n{ \ep^{-1} \nabla r_{\ep} }_{L^{\infty}}
  + \n{ \ep^{-1} \nabla q_{\ep} }_{L^{\infty}}
   }
   \sp{
   \n{\ep^{-1} \nabla^{s} R_{\ep}  }_{L^2}
   +
    \n{\ep^{-1} \nabla^{s} Q_{\ep}  }_{L^2}
   } \\
   & \quad 
   +
   C \n{\vu}_{H^s} 
   \sp{
  \n{ \ep^{-1} \nabla R_{\ep} }_{L^{\infty}}
  + \n{ \ep^{-1} \nabla Q_{\ep} }_{L^{\infty}}
   }  \\
   & \qquad  \qquad 
   \times 
   \sp{
   \n{ \ep^{-1} \nabla r_{\ep} }_{H^{s-1}}
   +
   \n{ \ep^{-1} \nabla q_{\ep} }_{H^{s-1}}
   + \ep^{s-1}  \n{ \ep^{-1} \nabla r_{\ep} }_{H^{s-1}}^s 
   +
   \ep^{s-1}  \n{ \ep^{-1} \nabla q_{\ep} }_{H^{s-1}}^s 
   } \\
   & \leq 
   C \sp{ 
   C_1^{\f{1}{2}} + \ep^{s-1} C_1^{\f{s}{2}}
   }
   \sp{
\n{ \ep^{-1} (R_{\ep}-1)   }_{H^s}^2 
+ \n{ \ep^{-1} (Q_{\ep}-1)   }_{H^s}^2  
+\n{  \vu }_{H^s}^2
}. 
\end{align} }}
Similarly,{\small{
\begin{align}
|I_7| & \leq 
C \n{\vu}_{H^s} \left[
\n{  \nabla \sp{ \f{1}{r_{\ep}+q_{\ep}} }  }_{L^{\infty}}
\n{ \nabla^{s-1}\Delta \vu }_{L^2}
+
\n{\Delta \vu }_{L^{\infty}} \n{ \nabla^{s} \sp{ \f{1}{r_{\ep}+q_{\ep}} }  }_{L^2}
\right] \\
& \quad 
+ 
C \n{\vu}_{H^s} \left[
\n{  \nabla \sp{ \f{1}{r_{\ep}+q_{\ep}} }  }_{L^{\infty}}
\n{ \nabla^{s-1} \nabla \div \vu }_{L^2}
+
\n{ \nabla \div \vu }_{L^{\infty}} \n{ \nabla^{s} \sp{ \f{1}{r_{\ep}+q_{\ep}} }  }_{L^2}
\right] \\
& 
\leq 
C \n{\vu}_{H^s}
\sp{   
\n{\nabla \vu}_{H^s} +\n{\div \vu}_{H^s}
}
\Big[
\sp{
\n{\nabla r_{\ep}}_{L^{\infty}}+ \n{\nabla q_{\ep}}_{L^{\infty}}
} \\
& \qquad \qquad
+  \sp{
\n{\nabla r_{\ep}}_{H^{s-1}}+ \n{\nabla q_{\ep}}_{H^{s-1}}
+\n{\nabla r_{\ep}}_{H^{s-1}}^{s} +\n{\nabla q_{\ep}}_{H^{s-1}}^{s}
}
\Big] \\
& 
\leq 
C \n{\vu}_{H^s}
\sp{   
\n{\nabla \vu}_{H^s} +\n{\div \vu}_{H^s}
}
\Big[
\ep
\sp{
 \n{ \ep^{-1} \nabla r_{\ep}}_{L^{\infty}}+ \n{\ep^{-1} \nabla q_{\ep}}_{L^{\infty}}
} \\
& \qquad \qquad
+  \ep \sp{
\n{ \ep^{-1} \nabla r_{\ep}}_{H^{s-1}}+ \n{ \ep^{-1} \nabla q_{\ep}}_{H^{s-1}}
}
+ \ep^{s}
\sp{
\n{ \ep^{-1} \nabla r_{\ep}}_{H^{s-1}}^{s} +\n{\ep^{-1} \nabla q_{\ep}}_{H^{s-1}}^{s}
}
\Big] \\
& \leq 
C \n{\vu}_{H^s}
\sp{   
\n{\nabla \vu}_{H^s} +\n{\div \vu}_{H^s}
}
\sp{
\ep C_1^{\f{1}{2}} + \ep^{s} C_1^{\f{s}{2}}
} \\
   & \leq 
   \f{\mu}{100} \n{\nabla \vu}_{H^s}^2 + 
   \f{\mu+\lambda}{100} \n{\div \vu}_{H^s}^2 
   + 
   C \n{\vu}_{H^s}^{2} \sp{
   \ep^2 C_1 +\ep^{2 s} C_1^{s}
   }.
\end{align}  }}

\noindent 
{\it Step 4. The estimates for $I_4-I_7$ when $D^{\alpha}=\p_t \nabla^{\beta},|\beta|\leq s-1$.} 
By \eqref{unif-bd}, Sobolev's embedding inequality, Lemma \ref{lemm:prod} and Lemma \ref{lemm:compo}, we see {\small{
\begin{align}
 |I_4| & \leq 
 C \n{\ep^{-1} \p_t  R_{\ep} }_{ H^{s-1} }
 \sp{
 \n{\p_t \vv}_{L^{\infty}} \n{ \ep^{-1} \nabla^{s-1} \nabla R_{\ep}  }_{L^2}
 +
 \n{ \ep^{-1} \nabla R_{\ep} }_{L^{\infty}}
 \n{ \nabla^{s-1} \p_t \vv }_{L^2}
 } \\
 & \quad 
 + 
 C \n{\ep^{-1} \p_t  R_{\ep} }_{ H^{s-1} }
 \sp{
 \n{\nabla \vv}_{L^{\infty}} 
 \n{\ep^{-1} \nabla^{s-2} \nabla \p_t R_{\ep} }_{L^2}
 + 
 \n{ \ep^{-1} \nabla  \p_t R_{\ep}  }_{L^{\infty}}
 \n{ \nabla^{s-1} \vv }_{L^2}
 } \\
 & \quad 
 + 
 C \n{\ep^{-1} \p_t  R_{\ep} }_{ H^{s-1} } 
 \sp{
 \n{ \ep^{-1}  \p_t r_{\ep}  }_{L^{\infty}}
 \n{ \nabla^{s-1} \div \vu  }_{L^2}
 +
 \n{\div \vu}_{L^{\infty}}
 \n{ \ep^{-1} \nabla^{s-1} \p_t r_{\ep}   }_{L^2}
 }\\
 & \quad 
 + 
 C \n{\ep^{-1} \p_t  R_{\ep} }_{ H^{s-1} } 
 \sp{
 \n{\ep^{-1} \nabla r_{\ep}  }_{L^{\infty}}
 \n{ \nabla^{s-2} \div \p_t \vu }_{L^2}
 +
 \n{  \div \p_t \vu }_{L^{\infty}}
 \n{ 
 \ep^{-1} \nabla^{s-1} r_{\ep}
 }_{L^2} 
 } \\ 
 & \quad 
 + 
 C \n{\ep^{-1} \p_t  Q_{\ep} }_{ H^{s-1} }
 \sp{
 \n{\p_t \vv}_{L^{\infty}} \n{ \ep^{-1} \nabla^{s-1} \nabla Q_{\ep}  }_{L^2}
 +
 \n{ \ep^{-1} \nabla Q_{\ep} }_{L^{\infty}}
 \n{ \nabla^{s-1} \p_t \vv }_{L^2}
 } \\
 & \quad 
 + 
 C \n{\ep^{-1} \p_t  Q_{\ep} }_{ H^{s-1} }
 \sp{
 \n{\nabla \vv}_{L^{\infty}} 
 \n{\ep^{-1} \nabla^{s-2} \nabla \p_t Q_{\ep} }_{L^2}
 + 
 \n{ \ep^{-1} \nabla  \p_t Q_{\ep}  }_{L^{\infty}}
 \n{ \nabla^{s-1} \vv }_{L^2}
 } \\ 
 & \quad 
 + 
 C \n{\ep^{-1} \p_t  Q_{\ep} }_{ H^{s-1} } 
 \sp{
 \n{ \ep^{-1}  \p_t q_{\ep}  }_{L^{\infty}}
 \n{ \nabla^{s-1} \div \vu  }_{L^2}
 +
 \n{\div \vu}_{L^{\infty}}
 \n{ \ep^{-1} \nabla^{s-1} \p_t q_{\ep}   }_{L^2}
 }\\
 & \quad  
 + 
 C \n{\ep^{-1} \p_t  Q_{\ep} }_{ H^{s-1} } 
 \sp{
 \n{\ep^{-1} \nabla q_{\ep}  }_{L^{\infty}}
 \n{ \nabla^{s-2} \div \p_t \vu }_{L^2}
 +
 \n{  \div \p_t \vu }_{L^{\infty}}
 \n{ 
 \ep^{-1} \nabla^{s-1} q_{\ep}
 }_{L^2} 
 } \\ 
 &  
 \leq 
 C \n{\ep^{-1} \p_t  R_{\ep} }_{ H^{s-1} }
 \Big(
 C_1^{\f{1}{2}} \n{ \ep^{-1} \p_t R_{\ep} }_{ H^{s-1} }
 + 
 C_1^{\f{1}{2}} \n{  \p_t \vu}_{ H^{s-1} } \\
 & \qquad \qquad \qquad \qquad \qquad \qquad
 + C_2^{\f{1}{2}} \n{ \ep^{-1} (R_{\ep}-1)  }_{ H^{s} }
 + C_2^{\f{1}{2}} \n{ \vu  }_{ H^{s} }
 \Big) \\
 & \quad 
 + 
  C \n{\ep^{-1} \p_t  Q_{\ep} }_{ H^{s-1} }
 \Big(
 C_1^{\f{1}{2}} \n{ \ep^{-1} \p_t Q_{\ep} }_{ H^{s-1} }
 + 
 C_1^{\f{1}{2}} \n{  \p_t \vu}_{ H^{s-1} } \\
 & \qquad \qquad \qquad \qquad \qquad \qquad
 + C_2^{\f{1}{2}} \n{ \ep^{-1} (Q_{\ep}-1)  }_{ H^{s} }
 + C_2^{\f{1}{2}} \n{ \vu  }_{ H^{s} }
 \Big)  \\
 & \leq 
 C \sp{1 +C_1+C_2  } 
 \sp{
 \n{\ep^{-1} \p_t  R_{\ep} }_{ H^{s-1} }^2 
 + 
 \n{\ep^{-1} \p_t  Q_{\ep} }_{ H^{s-1} }^2
 +
 \n{  \p_t \vu}_{ H^{s-1} }^2
 } \\
 & \qquad 
 + C \sp{
 \n{ \ep^{-1} (R_{\ep}-1)  }_{ H^{s} }^2 
 + 
  \n{ \ep^{-1} (Q_{\ep}-1)  }_{ H^{s} }^2
  + 
  \n{ \vu  }_{ H^{s} }^2
 }. 
\end{align} }}
Using similar ideas, we estimate $I_5$ as 
\begin{align}
|I_5| & \leq 
C \n{ \p_t \vu}_{ H^{s-1} } 
\sp{
\n{\p_t \vv }_{L^{\infty}}
\n{ \nabla^{s-1} \nabla \vu  }_{L^2}
+ 
\n{\nabla \vu}_{L^{\infty}}
\n{ \nabla^{s-1} \p_t \vv  }_{L^2}
} \\
& \quad 
+ 
C \n{ \p_t \vu}_{ H^{s-1} } 
\sp{
\n{\nabla \vv}_{L^{\infty}} 
\n{ \nabla^{s-2}  \nabla \p_t \vu  }_{L^2}
+
\n{\nabla \p_t \vu  }_{L^{\infty}}
\n{\nabla^{s-1}  \vv }_{L^2}
} \\
& \leq 
C \n{ \p_t \vu}_{ H^{s-1} } 
\sp{
C_1^{\f{1}{2}} \n{ \p_t \vu}_{ H^{s-1} }
+ 
C_2^{\f{1}{2}} \n{\vu}_{H^s}
} \\
& \leq 
C \sp{1 +C_1+C_2  } 
\n{ \p_t \vu}_{ H^{s-1} }^2
+ C \n{\vu}_{H^s}^2. 
\end{align}

Before turning to the estimate of $I_6$, we make a preliminary computation. Based on Lemma \ref{lemm:prod} and Lemma \ref{lemm:compo}, it follows that
\begin{align}\label{aux-cal}
& \n{ \nabla^{s-1} \p_t \sp{ \f{1}{r_{\ep}+q_{\ep}} } }_{L^2} \\
& \quad 
= 
\n{
\nabla^{s-1} \sp{ \f{\p_t r_{\ep}}{(r_{\ep}+q_{\ep})^2}  }  
+
\nabla^{s-1} \sp{ \f{\p_t q_{\ep}}{(r_{\ep}+q_{\ep})^2}  }  
}_{L^2} \\
& \quad 
\leq C 
\sp{
\n{\p_t r_{\ep}}_{L^{\infty}} 
\n{\nabla^{s-1} \sp{ \f{1}{(r_{\ep}+q_{\ep})^2}  }  }_{L^2}
+
\n{ \f{1}{(r_{\ep}+q_{\ep})^2}  }_{L^{\infty}}
\n{ \nabla^{s-1}  \p_t r_{\ep}  }_{L^2}
} \\
& \qquad 
+ C
\sp{
\n{\p_t q_{\ep}}_{L^{\infty}} 
\n{\nabla^{s-1} \sp{ \f{1}{(r_{\ep}+q_{\ep})^2}  }  }_{L^2}
+
\n{ \f{1}{(r_{\ep}+q_{\ep})^2}  }_{L^{\infty}}
\n{ \nabla^{s-1}  \p_t q_{\ep}  }_{L^2}
} \\
& \quad 
\leq C  
\sp{  
\n{ \p_t r_{\ep} }_{H^{s-1}}
+\n{ \p_t q_{\ep} }_{H^{s-1}}
}
+ 
C  
\sp{  
\n{ \p_t r_{\ep} }_{H^{s-1}}
+\n{ \p_t q_{\ep} }_{H^{s-1}}
} \\
& \qquad \qquad \qquad 
\times \sp{
\n{ \nabla r_{\ep} }_{H^{s-2}}
+
\n{ \nabla q_{\ep} }_{H^{s-2}}
+
\n{ \nabla r_{\ep} }_{H^{s-2}}^{s-1}
+
\n{ \nabla q_{\ep} }_{H^{s-2}}^{s-1}
} \\
& \quad 
\leq 
C
\sp{
\ep C_2^{\f{1}{2}} + \ep C_1^{\f{1}{2}} \ep C_2^{\f{1}{2}}
+ 
\ep^{s-1} C_1^{\f{s-1}{2}}
\ep C_2^{\f{1}{2}}
}.
\end{align}
Now we estimate $I_6$ as{\small{
\begin{align}
|I_6| & 
\leq 
C  \ep^{-2}\n{ \p_t \vu }_{H^{s-1}}
\Big(
\n{  
\p_t \sp{ \f{ \p_{R}p(\Z)(r_{\ep},q_{\ep} )  }{r_{\ep}+q_{\ep}}  }
}_{L^{\infty}} 
\n{\nabla^{s-1} \nabla R_{\ep} }_{L^2} \\
& \qquad \qquad \qquad \qquad \qquad 
+
\n{ \nabla R_{\ep} }_{L^{\infty}} 
\n{  \nabla^{s-1} \p_t \sp{ \f{ \p_{R}p(\Z)(r_{\ep},q_{\ep} )  }{r_{\ep}+q_{\ep}}  } }_{L^2}
\Big) \\
& \quad 
+ 
C  \ep^{-2}\n{ \p_t \vu }_{H^{s-1}}
\Big( 
\n{
\nabla  \sp{ \f{ \p_{R}p(\Z)(r_{\ep},q_{\ep} )  }{r_{\ep}+q_{\ep}}  } 
}_{L^{\infty}} 
\n{ \nabla^{s-2} \nabla \p_t R_{\ep} }_{L^2} \\
& \qquad \qquad \qquad \qquad \qquad 
+
\n{  \nabla \p_t R_{\ep}   }_{L^{\infty}} 
\n{  \nabla^{s-1}  \sp{ \f{ \p_{R}p(\Z)(r_{\ep},q_{\ep} )  }{r_{\ep}+q_{\ep}}  } }_{L^2}  \Big)   \\ 
& \quad 
 + 
 C  \ep^{-2}\n{ \p_t \vu }_{H^{s-1}}
\Big(
\n{  
\p_t \sp{ \f{ \p_{Q}p(\Z)(r_{\ep},q_{\ep} )  }{r_{\ep}+q_{\ep}}  }
}_{L^{\infty}} 
\n{\nabla^{s-1} \nabla Q_{\ep} }_{L^2} \\
& \qquad \qquad \qquad \qquad \qquad 
+
\n{ \nabla Q_{\ep} }_{L^{\infty}} 
\n{  \nabla^{s-1} \p_t \sp{ \f{ \p_{Q}p(\Z)(r_{\ep},q_{\ep} )  }{r_{\ep}+q_{\ep}}  } }_{L^2}
\Big) \\
& \quad 
+ 
C  \ep^{-2}\n{ \p_t \vu }_{H^{s-1}}
\Big( 
\n{
\nabla  \sp{ \f{ \p_{Q}p(\Z)(r_{\ep},q_{\ep} )  }{r_{\ep}+q_{\ep}}  } 
}_{L^{\infty}} 
\n{ \nabla^{s-2} \nabla \p_t Q_{\ep} }_{L^2} \\
& \qquad \qquad \qquad \qquad \qquad 
+
\n{  \nabla \p_t Q_{\ep}   }_{L^{\infty}} 
\n{  \nabla^{s-1}  \sp{ \f{ \p_{Q}p(\Z)(r_{\ep},q_{\ep} )  }{r_{\ep}+q_{\ep}}  } }_{L^2}  \Big)   \\ 
& 
\leq  
C \n{ \p_t \vu }_{H^{s-1}}
\sp{
\n{\ep^{-1} \nabla r_{\ep} }_{L^{\infty}} 
+
\n{\ep^{-1} \nabla q_{\ep} }_{L^{\infty}} 
}
\sp{
\n{\ep^{-1} \nabla R_{\ep} }_{H^{s-1}}
+
\n{\ep^{-1} \nabla Q_{\ep} }_{H^{s-1}}
} \\
& \quad 
+ 
C \n{ \p_t \vu }_{H^{s-1}}
\sp{
\n{\ep^{-1} \p_t r_{\ep}}_{L^{\infty}} 
+
\n{\ep^{-1} \p_t q_{\ep}}_{L^{\infty}} 
}
\sp{
\n{\ep^{-1}(R_{\ep}-1)}_{H^{s}}
+ 
\n{\ep^{-1}(Q_{\ep}-1)}_{H^{s}}
} \\
& \quad 
+ 
C \n{ \p_t \vu }_{H^{s-1}} 
\sp{
\n{ \ep^{-1} \p_t R_{\ep}}_{H^{s-1}} 
+ \n{ \ep^{-1} \p_t Q_{\ep}}_{H^{s-1}} 
} \\ 
& \qquad \quad 
\times 
\sp{
\n{ \ep^{-1} \p_t r_{\ep}}_{H^{s-2}} 
+\n{ \ep^{-1} \p_t q_{\ep}}_{H^{s-2}} 
+ \ep^{s-2} \n{ \ep^{-1} \p_t r_{\ep}}_{H^{s-2}}^{s-1}
+ \ep^{s-2} \n{ \ep^{-1} \p_t q_{\ep}}_{H^{s-2}}^{s-1}
} \\
& \quad 
+ C \n{ \p_t \vu }_{H^{s-1}}  
\n{\ep^{-1}(R_{\ep}-1)}_{H^{s}} 
\n{  \nabla^{s-1} \p_t \sp{ \f{ \p_{R}p(\Z)(r_{\ep},q_{\ep} )  }{r_{\ep}+q_{\ep}}  } }_{L^2} \\
&\quad 
+ C \n{ \p_t \vu }_{H^{s-1}}  
\n{\ep^{-1}(Q_{\ep}-1)}_{H^{s}} 
\n{  \nabla^{s-1} \p_t \sp{ \f{ \p_{Q}p(\Z)(r_{\ep},q_{\ep} )  }{r_{\ep}+q_{\ep}}  } }_{L^2} \\
& \leq 
C \sp{
   1+C_1^s +C_2^s
   } 
   \sp{
   \n{\p_t \vu }_{H^{s-1}}^2 
   + \n{\ep^{-1} \p_t R_{\ep}  }_{H^{s-1}}^2 
   + \n{\ep^{-1} \p_t Q_{\ep}  }_{H^{s-1}}^2
   }  \\
   & \qquad 
   +
   C \sp{
   \n{\ep^{-1}(R_{\ep}-1)}_{H^{s}}^2
   +  \n{\ep^{-1}(Q_{\ep}-1)}_{H^{s}}^2
   }, 
\end{align}  }}
where in the last step we essentially used the same arguments as \eqref{aux-cal}.

In the same spirit, we estimate $I_7$ as {\small{
\begin{align}
|I_7| & \leq 
C \n{ \p_t \vu}_{ H^{s-1} }
\sp{
\n{ \p_t \sp{ \f{1}{r_{\ep}+q_{\ep}}} }_{L^{\infty}}
\n{ \nabla^{s-1} \Delta \vu}_{L^2}
+
\n{\Delta \vu}_{L^{\infty}} 
\n{ \nabla^{s-1}  \p_t \sp{ \f{1}{r_{\ep}+q_{\ep}}} }_{L^2}
} \\
& \quad 
+ 
C \n{ \p_t \vu}_{ H^{s-1} }
\sp{
\n{ \nabla \sp{ \f{1}{r_{\ep}+q_{\ep}}} }_{L^{\infty}}
\n{  \nabla^{s-2}  \Delta \p_t \vu }_{L^2}
+ \n{  \Delta \p_t \vu }_{L^{\infty}}
\n{ \nabla^{s-1}   \sp{ \f{1}{r_{\ep}+q_{\ep}}} }_{L^2}
} \\ 
& 
\quad +
C \n{ \p_t \vu}_{ H^{s-1} }
\sp{
\n{ \p_t \sp{ \f{1}{r_{\ep}+q_{\ep}}} }_{L^{\infty}}
\n{ \nabla^{s-1} \nabla \div \vu}_{L^2}
+
\n{\nabla \div \vu}_{L^{\infty}} 
\n{ \nabla^{s-1}  \p_t \sp{ \f{1}{r_{\ep}+q_{\ep}}} }_{L^2}
} \\
& \quad 
+ 
C \n{ \p_t \vu}_{ H^{s-1} }
\sp{
\n{ \nabla \sp{ \f{1}{r_{\ep}+q_{\ep}}} }_{L^{\infty}}
\n{  \nabla^{s-2}  \nabla \div \p_t \vu }_{L^2}
+ \n{  \nabla \div \p_t \vu }_{L^{\infty}}
\n{ \nabla^{s-1}   \sp{ \f{1}{r_{\ep}+q_{\ep}}} }_{L^2}
} \\ 
& 
\leq 
C \n{ \p_t \vu}_{ H^{s-1} }
\sp{
\ep C_1^{\f{1}{2}} +\ep^{s-1} C_1^{\f{s-1}{2}}
}
\sp{
\n{\nabla \p_t \vu}_{ H^{s-1} }
+ \n{\div \p_t \vu}_{ H^{s-1} }
} \\
& \quad 
+ 
C \n{ \p_t \vu}_{ H^{s-1} } 
\sp{
\ep C_2^{\f{1}{2}} + \ep C_1^{\f{1}{2}} \ep C_2^{\f{1}{2}}
+ 
\ep^{s-1} C_1^{\f{s-1}{2}}
\ep C_2^{\f{1}{2}}
}
\sp{
\n{\nabla \vu}_{H^s} + \n{\div \vu}_{H^s}
} \\
& \leq  
\f{\mu}{100} \n{\nabla \p_t \vu}_{H^{s-1}}^2 + 
   \f{\mu+\lambda}{100} \n{\div \p_t \vu}_{H^{s-1}}^2 
   + 
   C \sp{
   1+C_1^s +C_2^s
   }\n{ \p_t \vu}_{ H^{s-1} }^2 \\
   & \qquad 
   +
   C \sp{
\n{\nabla \vu}_{H^s}^2 + \n{\div \vu}_{H^s}^2
}.
\end{align} }}

\subsection{Closure of the energy estimates}\label{subs:clo}
Combining the estimates in Step 1 and Step 3, we obtain{\small{
\begin{align}\label{eq:ene-est-1}
& \f{1}{2} \f{d}{dt} \sum_{|\alpha|\leq s}
\int_{\Om} 
\Big[
\frac{ \gamma_{+} [ \Z^{\gamma_{+}-1} (\p_{R}\Z ) ](r_{\ep},q_{\ep})    }{\ep^2   r_{\ep}  } | \nabla^{\alpha}(R_{\ep}-1)|^2 \\
&  \qquad \qquad +
\frac{ \gamma_{+} [ \Z^{\gamma_{+}-1} (\p_{Q}\Z ) ](r_{\ep},q_{\ep})    }{\ep^2    q_{\ep} } |\nabla^{\alpha}(Q_{\ep}-1)|^2
+ 
(r_{\ep}+q_{\ep}) |\nabla^{\alpha} \vu|^2
\Big] dx \\
& \qquad \qquad +
\mu \sum_{|\alpha|\leq s} \n{ \Grad \nabla^{\alpha}\vu  }_{L^2}^2 
+ (\mu+\lambda) \sum_{|\alpha|\leq s} \n{ \div \nabla^{\alpha}\vu  }_{L^2}^2 \\
& \quad 
\leq 
C \sp{
C_1
+
C_1^{\f{1}{2}}
+
C_1^{\f{s}{2}}
+
C_1^{s}
+
C_2^{\f{1}{2}}
}
\sp{
\n{ \ep^{-1}(R_{\ep}-1) }_{H^s}^2
+ 
\n{ \ep^{-1}(Q_{\ep}-1) }_{H^s}^2
+ 
\n{ \vu}_{H^s}^2
}
.
\end{align}  }}
From the assumptions of Theorem \ref{Thm1} we see that
\begin{align}
\mathscr{E}_{s} (R_{0,\ep},Q_{0,\ep},\mathbf{u}_{0,\ep})
&
= \f{1}{2}
\sp{
\n{ \ep^{-1}  \widetilde{R_{0,\ep}} }_{H^s}^2
+
\n{ \ep^{-1}  \widetilde{Q_{0,\ep}} }_{H^s}^2
+
\n{ \mathbf{u}_0+ \widetilde{\mathbf{u}_{0,\ep}} }_{H^s}^2
} \\
& \leq 
C
\sp{
\n{\mathbf{u}_0  }_{H^s}^2
+ \ep^2 \delta_0^2
}.
\end{align}
By the uniform bounds \eqref{unif-bd} and Gr\"{o}nwall's inequality, we infer for all $t\in [0,T_{\ast}]$ that 
\begin{align}
& \mathscr{E}_{s} (R_{\ep},Q_{\ep},\vu)(t)
+ 
\mu \int_0^t \n{\nabla \vu}_{H^s}^2 d\tau
+
(\mu+\lambda) \int_0^t \n{\div \vu}_{H^s}^2 d\tau \\
& \quad 
\leq C \exp \left[  
C \sp{
C_1
+
C_1^{\f{1}{2}}
+
C_1^{\f{s}{2}}
+
C_1^{s}
+
C_2^{\f{1}{2}}
} t
\right] 
\sp{
\n{\mathbf{u}_0  }_{H^s}^2
+ \ep^2 \delta_0^2
}. 
\end{align}
Hence, if $T_{\ast}$ is sufficiently small such that
\begin{align}\label{cond-T}
  T_{\ast}
  \leq \f{1}{
  C_1
+
C_1^{\f{1}{2}}
+
C_1^{\f{s}{2}}
+
C_1^{s}
+
C_2^{\f{1}{2}}
  },
\end{align}
then we conclude  for all $t\in [0,T_{\ast}]$ that
\begin{align}\label{eq:clo-1}
& \mathscr{E}_{s} (R_{\ep},Q_{\ep},\vu)(t)
+ 
\mu \int_0^t \n{\nabla \vu}_{H^s}^2 d\tau
+
(\mu+\lambda) \int_0^t \n{\div \vu}_{H^s}^2 d\tau \\
& \quad 
\leq C
\sp{
\n{\mathbf{u}_0  }_{H^s}^2
+ \ep^2 \delta_0^2
} 
\leq 
C \sp{ 1+ \n{\mathbf{u}_0  }_{H^s}^2 }
\equiv C_1,
\end{align}
as long as $\ep$ and $\delta_0$ are small enough.

To proceed, we combine Step 1, Step 2 and Step 4 to deduce that
\begin{align}
 & \f{1}{2} \f{d}{dt} \sum_{|\beta|\leq s-1}
\int_{\Om} 
\Big[
\frac{ \gamma_{+} [ \Z^{\gamma_{+}-1} (\p_{R}\Z ) ](r_{\ep},q_{\ep})    }{\ep^2   r_{\ep}  } | \nabla^{\beta} \p_t R_{\ep}|^2 \\
&  \qquad \qquad +
\frac{ \gamma_{+} [ \Z^{\gamma_{+}-1} (\p_{Q}\Z ) ](r_{\ep},q_{\ep})    }{\ep^2    q_{\ep} } | \nabla^{\beta} \p_tQ_{\ep}|^2
+ 
(r_{\ep}+q_{\ep}) | \nabla^{\beta} \p_t \vu|^2
\Big] dx \\
& \qquad \qquad +
\mu \sum_{|\beta|\leq s-1} \n{ \nabla^{\beta} \p_t\vu  }_{L^2}^2 
+ (\mu+\lambda) \sum_{|\beta|\leq s-1} \n{ \div 
\nabla^{\beta} \p_t\vu  }_{L^2}^2 \\ 
& \quad 
\leq 
C \sp{
1+C_1^s + C_2^s
}
\sp{
   \n{\p_t \vu }_{H^{s-1}}^2 
   + \n{\ep^{-1} \p_t R_{\ep}  }_{H^{s-1}}^2 
   + \n{\ep^{-1} \p_t Q_{\ep}  }_{H^{s-1}}^2
   } \\
   & \qquad 
   + 
   C \sp{
   1+ \n{\mathbf{u}_0  }_{H^s}^2 +
   \n{\nabla \vu}_{H^s}^2 + \n{\div \vu}_{H^s}^2
   }, 
\end{align}
where we have used \eqref{eq:clo-1} in the last step. 
Again from the assumptions of Theorem \ref{Thm1}, the governing equations \eqref{eq:two-fluid-1} and the uniform bounds \eqref{unif-bd}, we see that
\begin{align}
& \mathscr{E}_{s-1} (\p_t R_{\ep},\p_t Q_{\ep},\p_t\vu)(0)   \\
& \quad 
\leq 
C 
\Big(   
\n{
\ep^{-1} \sp{ \mathbf{u}_0 + \widetilde{\mathbf{u}_{0,\ep}}  }
\cdot \nabla \widetilde{R_{0,\ep}}
}_{H^{s-1}}^2
+
\n{
\ep^{-1} \sp{ \mathbf{u}_0 + \widetilde{\mathbf{u}_{0,\ep}}  }
\cdot \nabla \widetilde{Q_{0,\ep}}
}_{H^{s-1}}^2 \\
& \qquad 
+ 
\n{
\ep^{-1} \sp{ 1+\widetilde{R_{0,\ep}}  } \div \widetilde{\mathbf{u}_{0,\ep}}
}_{H^{s-1}}^2
+ 
\n{
\ep^{-1} \sp{ 1+\widetilde{Q_{0,\ep}}  } \div \widetilde{\mathbf{u}_{0,\ep}}
}_{H^{s-1}}^2 \\
& \qquad 
+ 
\n{
\sp{ \mathbf{u}_0 + \widetilde{\mathbf{u}_{0,\ep}}  } 
\cdot 
\nabla
\sp{ \mathbf{u}_0 + \widetilde{\mathbf{u}_{0,\ep}}  }
}_{H^{s-1}}^2  
+
\n{
\ep^{-2} \nabla \widetilde{R_{0,\ep}} 
}_{H^{s-1}}^2 
+
\n{
\ep^{-2} \nabla \widetilde{Q_{0,\ep}} 
}_{H^{s-1}}^2  \\
& \qquad 
+   
\n{
\Delta \sp{ \mathbf{u}_0 + \widetilde{\mathbf{u}_{0,\ep}}  } 
+\nabla \div \sp{ \mathbf{u}_0 + \widetilde{\mathbf{u}_{0,\ep}}  } 
}_{H^{s-1}}^2
\Big)    \\
& \quad 
\leq 
C \sp{ 
\n{\mathbf{u}_0  }_{H^{s+1}}^2 
+\delta_0^2 
}, 
\end{align}
from which and the uniform bounds \eqref{unif-bd} and Gr\"{o}nwall's inequality, we infer for all $t\in [0,T_{\ast}]$ that 
\begin{align}
& \mathscr{E}_{s-1} (\p_t R_{\ep},\p_t Q_{\ep},\p_t\vu)(t)
+ 
\mu \int_0^t \n{\nabla \p_t\vu}_{H^{s-1}}^2 d\tau
+
(\mu+\lambda) \int_0^t \n{\div \p_t \vu}_{H^{s-1}}^2 d\tau \\
& \quad     
\leq 
C \exp \left[  
C \sp{
1+C_1^s + C_2^s
} t 
\right] \\
& \qquad \qquad  
\times 
\sp{
\n{\mathbf{u}_0  }_{H^{s+1}}^2 
+\delta_0^2 
+ 
\int_0^t 
\sp{
   1+ \n{\mathbf{u}_0  }_{H^s}^2 +
   \n{\nabla \vu}_{H^s}^2 + \n{\div \vu}_{H^s}^2
   }
   d\tau
}
\end{align}
Therefore, if $T_{\ast}$ is small enough such that
\begin{align}
  T_{\ast}
  \leq 
  \min \left\{ 
  \f{1}{
  C_1
+
C_1^{\f{1}{2}}
+
C_1^{\f{s}{2}}
+
C_1^{s}
+
C_2^{\f{1}{2}}
  },
  \f{1}{1+C_1^s + C_2^s},
  1
  \right\}
  ,
\end{align}
then we conclude for all $t\in [0,T_{\ast}]$ that 
\begin{align}
 & \mathscr{E}_{s-1} (\p_t R_{\ep},\p_tQ_{\ep},\p_t\vu)(t)
+ 
\mu \int_0^t \n{\nabla \p_t\vu}_{H^{s-1}}^2 d\tau
+
(\mu+\lambda) \int_0^t \n{\div \p_t \vu}_{H^{s-1}}^2 d\tau 
\nonumber \\
& \quad        
\leq 
C \sp{
\n{\mathbf{u}_0  }_{H^{s+1}}^2 
+\delta_0^2 
+
1+ \n{\mathbf{u}_0  }_{H^s}^2
}
\leq C \sp{ 1+ \n{\mathbf{u}_0  }_{H^{s+1}}^2 }
\equiv C_2 ,   \label{eq:clo-2}
\end{align}
as long as $\ep$ and $\delta_0$ are small enough.

Finally, we show that $R_{\ep},Q_{\ep}$ are close to $1$. To this end, we see from \eqref{eq:lin} that 
{\small{
\begin{align}\label{eq:lin-2}
    \begin{dcases}
        \partial_t (R_{\ep}-1) + \vv \cdot \Grad (R_{\ep}-1)+ r_{\ep} \div (\vu- \mathbf{u}_0) = 0, \\
         \partial_t  (Q_{\ep}-1) + \vv \cdot \Grad (Q_{\ep}-1) + q_{\ep} \div (\vu- \mathbf{u}_0) = 0, \\
        \p_t (\vu- \mathbf{u}_0) + \vv \cdot \Grad (\vu- \mathbf{u}_0)  + 
        \frac{ \gamma_{+} [ \Z^{\gamma_{+}-1} (\p_{R}\Z ) ](r_{\ep},q_{\ep})    }{\ep^2   (r_{\ep} + q_{\ep}) }
        \Grad ( R_{\ep}  -1) \\
        \quad 
        +
        \frac{ \gamma_{+} [ \Z^{\gamma_{+}-1} (\p_{Q}\Z ) ](r_{\ep},q_{\ep})    }{\ep^2   (r_{\ep} + q_{\ep}) }
        \Grad (Q_{\ep} -1)
        = \f{\mu}{  r_{\ep} + q_{\ep} } \Delta (\vu- \mathbf{u}_0)     
        +\f{  \mu+\lambda  }{  r_{\ep} + q_{\ep} } \nabla \div (\vu- \mathbf{u}_0) \\
        \quad+
        \f{\mu}{  r_{\ep} + q_{\ep} } \Delta \mathbf{u}_0
        - \vv \cdot \nabla \mathbf{u}_0. 
    \end{dcases}
\end{align}  }}
Following step by step the derivation of \eqref{eq:ene-est-1}, one obtains
{\small{
\begin{align}\label{eq:ene-est-2}
& \f{1}{2} \f{d}{dt} \sum_{|\alpha|\leq s}
\int_{\Om} 
\Big[
\frac{ \gamma_{+} [ \Z^{\gamma_{+}-1} (\p_{R}\Z ) ](r_{\ep},q_{\ep})    }{\ep^2   r_{\ep}  } | \nabla^{\alpha}(R_{\ep}-1)|^2 \\
&  \qquad \qquad +
\frac{ \gamma_{+} [ \Z^{\gamma_{+}-1} (\p_{Q}\Z ) ](r_{\ep},q_{\ep})    }{\ep^2    q_{\ep} } |\nabla^{\alpha}(Q_{\ep}-1)|^2
+ 
(r_{\ep}+q_{\ep}) |\nabla^{\alpha} (\vu-\mathbf{u}_0 )|^2
\Big] dx \\
& \qquad \qquad +
\mu \sum_{|\alpha|\leq s} \n{ \Grad \nabla^{\alpha} (\vu-\mathbf{u}_0 )  }_{L^2}^2 
+ (\mu+\lambda) \sum_{|\alpha|\leq s} \n{ \div \nabla^{\alpha} (\vu -\mathbf{u}_0 ) }_{L^2}^2 \\
& \quad 
\leq 
C \sp{
C_1
+
C_1^{\f{1}{2}}
+
C_1^{\f{s}{2}}
+
C_1^{s}
+
C_2^{\f{1}{2}}
}  \\
& \qquad \qquad \times 
\sp{
\n{ \ep^{-1}(R_{\ep}-1) }_{H^s}^2
+ 
\n{ \ep^{-1}(Q_{\ep}-1) }_{H^s}^2
+ 
\n{ (\vu-\mathbf{u}_0 )}_{H^s}^2
} \\
& \qquad \qquad
+
C \sp{1+ C_1^s } \n{\mathbf{u}_0}_{H^{s+1}}^2
.
\end{align}  }}
By the assumptions of Theorem \ref{Thm1} and Gr\"{o}nwall's inequality, we infer for all $t\in [0,T_{\ast}]$ that 
\begin{align}
& \mathscr{E}_{s} (R_{\ep},Q_{\ep},\vu-\mathbf{u}_0)(t)
+ 
\mu \int_0^t \n{\nabla (\vu-\mathbf{u}_0)}_{H^s}^2 d\tau
+
(\mu+\lambda) \int_0^t \n{\div (\vu-\mathbf{u}_0) }_{H^s}^2 d\tau \\
& \quad 
\leq C \exp \left[  
C \sp{
C_1
+
C_1^{\f{1}{2}}
+
C_1^{\f{s}{2}}
+
C_1^{s}
+
C_2^{\f{1}{2}}
} t
\right] 
\Big(
C \sp{ 1+C_1^s}t
+ \ep^2 \delta_0^2
\Big). 
\end{align}
Hence, in addition to the constraint \eqref{cond-T}, we may take $T_{\ast},\delta_0,\ep$ small enough such that 
\begin{align}
 & \mathscr{E}_{s} (R_{\ep},Q_{\ep},\vu-\mathbf{u}_0)(t) \\
 & \quad 
 \leq 
 C \exp \left[  
C \sp{
C_1
+
C_1^{\f{1}{2}}
+
C_1^{\f{s}{2}}
+
C_1^{s}
+
C_2^{\f{1}{2}}
} T_{\ast}
\right] 
\Big(
C \sp{ 1+C_1^s}T_{\ast}
+ \ep^2 \delta_0^2
\Big)
\leq
\delta, 
\end{align}
which together with Sobolev's embedding inequality gives
\begin{align}\label{eq:clo-3}
     | R_{\ep}-1|+|Q_{\ep}-1| < \ep \delta.
\end{align}

Combining \eqref{eq:clo-1}, \eqref{eq:clo-2} and \eqref{eq:clo-3}, we conclude that 
\begin{align}
    (R_{\ep}, Q_{\ep}, \vu)
    \equiv
    \chi (r_{\ep},q_{\ep},\vv)
    \in
    \mathcal{A}_{T_{\ast},\ep  }(R_{0,\ep}, Q_{0,\ep}, \mathbf{u}_{0,\ep},\delta,C_1,C_2).
\end{align}

\subsection{Contraction of the map \texorpdfstring{$\chi$}{}}\label{subs:cont}
The aim of this subsection is to prove that $\chi$ is indeed a contraction map in $\mathcal{A}_{T_{\ast},\ep  }(R_{0,\ep}, Q_{0,\ep}, \mathbf{u}_{0,\ep},\delta,C_1,C_2)$. More precisely, suppose that $(R_{\ep}, Q_{\ep}, \vu)= \chi (r_{\ep},q_{\ep},\vv)$ and $(R_{\ep}^{\#}, Q_{\ep}^{\#}, \vu^{\#})= \chi (r_{\ep}^{\#},q_{\ep}^{\#},\vv^{\#})$, where $(r_{\ep},q_{\ep},\vv),(r_{\ep}^{\#},q_{\ep}^{\#},\vv^{\#}) \in \mathcal{A}_{T_{\ast},\ep  }(R_{0,\ep}, Q_{0,\ep}, \mathbf{u}_{0,\ep},\delta,C_1,C_2)$. We shall prove that for $T_{\ast}$ small enough it holds
\begin{align}\label{eq:cont-map}
&\sup_{0 \leq t \leq T_{\ast}}
\sp{
\n{\ep^{-1}(R_{\ep}-R_{\ep}^{\#}) }_{L^2}^2
+
\n{\ep^{-1}(Q_{\ep}-Q_{\ep}^{\#}) }_{L^2}^2
+
\n{ \vu-\vu^{\#}}_{L^2}^2
} \\
& \qquad \qquad
+
\int_0^t \n{ \vu-\vu^{\#}}_{H^1}^2 d\tau  \\
& \quad \leq 
\vartheta 
\sup_{0 \leq t \leq T_{\ast}}
\sp{
\n{\ep^{-1}(r_{\ep}-r_{\ep}^{\#}) }_{L^2}^2
+
\n{\ep^{-1}(q_{\ep}-q_{\ep}^{\#}) }_{L^2}^2
+
\n{ \vv-\vv^{\#}}_{L^2}^2
}
\end{align}
for some $\vartheta \in (0,1)$. 

Direct calculation shows that $(R_{\ep}-R_{\ep}^{\#},Q_{\ep}-Q_{\ep}^{\#},\vu-\vu^{\#})$ solves
\begin{align}\label{eq:diffe}
    \begin{dcases}
        \partial_t  (R_{\ep}-R_{\ep}^{\#})
        + \vv \cdot \nabla (R_{\ep}-R_{\ep}^{\#})
        + (\vv-\vv^{\#}) \cdot \nabla R_{\ep}^{\#}
        +r_{\ep} \div (\vu-\vu^{\#}) \\
        \qquad +
        (r_{\ep}-r_{\ep}^{\#}) \div \vu^{\#}=0
        ,
       \\
       \partial_t  (Q_{\ep}-Q_{\ep}^{\#})
        + \vv \cdot \nabla (Q_{\ep}-Q_{\ep}^{\#})
        + (\vv-\vv^{\#}) \cdot \nabla Q_{\ep}^{\#}
        +q_{\ep} \div (\vu-\vu^{\#}) \\
        \qquad +
        (q_{\ep}-q_{\ep}^{\#}) \div \vu^{\#}=0
        ,
       \\ 
       \p_t (\vu-\vu^{\#}) 
       + \vv \cdot \nabla (\vu-\vu^{\#}) 
       +
        \frac{ \gamma_{+} [ \Z^{\gamma_{+}-1} (\p_{R}\Z ) ](r_{\ep},q_{\ep})    }{\ep^2   (r_{\ep} + q_{\ep}) }
        \Grad (R_{\ep}-R_{\ep}^{\#})  \\
        \qquad 
        +
        \frac{ \gamma_{+} [ \Z^{\gamma_{+}-1} (\p_{Q}\Z ) ](r_{\ep},q_{\ep})    }{\ep^2   (r_{\ep} + q_{\ep}) }
        \Grad (Q_{\ep}-Q_{\ep}^{\#}) 
        + 
        (\vv-\vv^{\#}) \cdot \nabla \vu^{\#}  \\
        \qquad
        + 
        \sp{
        \frac{ \gamma_{+} [ \Z^{\gamma_{+}-1} (\p_{R}\Z ) ](r_{\ep},q_{\ep})    }{\ep^2   (r_{\ep} + q_{\ep}) }
        - 
        \frac{ \gamma_{+} [ \Z^{\gamma_{+}-1} (\p_{R}\Z ) ](r_{\ep}^{\#},q_{\ep}^{\#})    }{\ep^2   (r_{\ep}^{\#} + q_{\ep}^{\#}) }
        }
        \Grad R_{\ep}^{\#}  \\
        \qquad
        + 
        \sp{
        \frac{ \gamma_{+} [ \Z^{\gamma_{+}-1} (\p_{Q}\Z ) ](r_{\ep},q_{\ep})    }{\ep^2   (r_{\ep} + q_{\ep}) }
        - 
        \frac{ \gamma_{+} [ \Z^{\gamma_{+}-1} (\p_{Q}\Z ) ](r_{\ep}^{\#},q_{\ep}^{\#})    }{\ep^2   (r_{\ep}^{\#} + q_{\ep}^{\#}) }
        }
        \Grad Q_{\ep}^{\#}  \\
        \quad 
        =
        \f{\mu}{  r_{\ep} + q_{\ep} } \Delta (\vu-\vu^{\#}) 
        +\f{\mu+\lambda}{  r_{\ep} + q_{\ep} }  \nabla \div (\vu-\vu^{\#}) \\
        \qquad 
        + 
        \sp{ 
         \f{\mu}{  r_{\ep} + q_{\ep} } -
          \f{\mu}{  r_{\ep}^{\#} + q_{\ep}^{\#} }
        } \Delta \vu^{\#}
        + 
        \sp{
        \f{\mu+\lambda}{  r_{\ep} + q_{\ep} }
        -
        \f{\mu+\lambda}{  r_{\ep}^{\#} + q_{\ep}^{\#} }
        } \nabla \div \vu^{\#}.
    \end{dcases}
\end{align}
To proceed, we take $L^2$-inner product of \eqref{eq:diffe}$_1$ with $\frac{ \gamma_{+} [ \Z^{\gamma_{+}-1} (\p_{R}\Z ) ](r_{\ep},q_{\ep})    }{\ep^2   r_{\ep}  } (R_{\ep}-R_{\ep}^{\#})$, \eqref{eq:diffe}$_2$ with $\frac{ \gamma_{+} [ \Z^{\gamma_{+}-1} (\p_{Q}\Z ) ](r_{\ep},q_{\ep})    }{\ep^2    q_{\ep} } (Q_{\ep}-Q_{\ep}^{\#})$ and\eqref{eq:diffe}$_3$ with $(r_{\ep}+q_{\ep})  (\vu-\vu^{\#}) $. Then, we add them together and employ integration by parts and some cancellations. By similar arguments as in Subsection \ref{subs:unif-est}, we arrive at {\small{
\begin{align}
& \f{1}{2} \f{d}{dt} 
\int_{\Om} 
\Big[
\frac{ \gamma_{+} [ \Z^{\gamma_{+}-1} (\p_{R}\Z ) ](r_{\ep},q_{\ep})    }{\ep^2   r_{\ep}  } (R_{\ep}-R_{\ep}^{\#})^2 \\
&  \qquad  +
\frac{ \gamma_{+} [ \Z^{\gamma_{+}-1} (\p_{Q}\Z ) ](r_{\ep},q_{\ep})    }{\ep^2    q_{\ep} } (Q_{\ep}-Q_{\ep}^{\#})^2
+ 
(r_{\ep}+q_{\ep}) (\vu-\vu^{\#})^2
\Big] dx \\
& \qquad  +
\mu \n{ \Grad (\vu-\vu^{\#})  }_{L^2}^2 
+ (\mu+\lambda) \n{ \div (\vu-\vu^{\#})  }_{L^2}^2 \\
& \quad 
= 
\f{1}{2} 
\int_{\Om}
\f{1}{\ep^2} (R_{\ep}-R_{\ep}^{\#})^2
\p_t \sp{
\frac{ \gamma_{+} [ \Z^{\gamma_{+}-1} (\p_{R}\Z ) ](r_{\ep},q_{\ep})    }{  r_{\ep}  }
} dx \\
& \qquad 
+
\f{1}{2} 
\int_{\Om}
\f{1}{\ep^2} (Q_{\ep}-Q_{\ep}^{\#})^2
\p_t \sp{
\frac{ \gamma_{+} [ \Z^{\gamma_{+}-1} (\p_{Q}\Z ) ](r_{\ep},q_{\ep})    }{  q_{\ep}  }
} dx \\
& \qquad 
+ 
\f{1}{2} 
\int_{\Om} (R_{\ep}-R_{\ep}^{\#})^2
\div \sp{ \frac{ \gamma_{+} [ \Z^{\gamma_{+}-1} (\p_{R}\Z ) ](r_{\ep},q_{\ep})    }{\ep^2   r_{\ep}  }  \vv   } dx \\
& \qquad 
+ 
\f{1}{2} 
\int_{\Om} (Q_{\ep}-Q_{\ep}^{\#})^2
\div \sp{ \frac{ \gamma_{+} [ \Z^{\gamma_{+}-1} (\p_{Q}\Z ) ](r_{\ep},q_{\ep})    }{\ep^2   q_{\ep}  }  \vv   } dx \\
& \qquad 
+
\f{1}{2} 
\int_{\Om}  | \vu-\vu^{\#}|^2 
\div ( (r_{\ep}+q_{\ep}) \vv  )  dx
+
\f{1}{2} 
\int_{\Om}  | \vu-\vu^{\#}|^2 
\p_t  (r_{\ep}+q_{\ep})   dx
\\
& \qquad
-
\int_{\Om}  
(r_{\ep}+q_{\ep}) 
\sp{
        \frac{ \gamma_{+} [ \Z^{\gamma_{+}-1} (\p_{R}\Z ) ](r_{\ep},q_{\ep})    }{\ep^2   (r_{\ep} + q_{\ep}) }
        - 
        \frac{ \gamma_{+} [ \Z^{\gamma_{+}-1} (\p_{R}\Z ) ](r_{\ep}^{\#},q_{\ep}^{\#})    }{\ep^2   (r_{\ep}^{\#} + q_{\ep}^{\#}) }
        }
        \Grad R_{\ep}^{\#} \cdot  (\vu-\vu^{\#}) dx \\
 & \qquad
-
\int_{\Om}  
(r_{\ep}+q_{\ep}) 
\sp{
        \frac{ \gamma_{+} [ \Z^{\gamma_{+}-1} (\p_{Q}\Z ) ](r_{\ep},q_{\ep})    }{\ep^2   (r_{\ep} + q_{\ep}) }
        - 
        \frac{ \gamma_{+} [ \Z^{\gamma_{+}-1} (\p_{Q}\Z ) ](r_{\ep}^{\#},q_{\ep}^{\#})    }{\ep^2   (r_{\ep}^{\#} + q_{\ep}^{\#}) }
        }
        \Grad Q_{\ep}^{\#} \cdot  (\vu-\vu^{\#}) dx \\   
  & \qquad
  +  \int_{\Om}  (r_{\ep}+q_{\ep}) 
  \sp{ 
         \f{\mu}{  r_{\ep} + q_{\ep} } -
          \f{\mu}{  r_{\ep}^{\#} + q_{\ep}^{\#} }
        } \Delta \vu^{\#} \cdot  (\vu-\vu^{\#}) dx \\   
  & \qquad
  + \int_{\Om}  (r_{\ep}+q_{\ep}) 
        \sp{
        \f{\mu+\lambda}{  r_{\ep} + q_{\ep} }
        -
        \f{\mu+\lambda}{  r_{\ep}^{\#} + q_{\ep}^{\#} }
        } \nabla \div \vu^{\#}  \cdot  (\vu-\vu^{\#}) dx \\
  & \qquad   
      - 
      \int_{\Om}  (r_{\ep}+q_{\ep}) 
      (\vv-\vv^{\#}) \cdot \nabla \vu^{\#}   \cdot  (\vu-\vu^{\#}) dx 
        \\
  & \qquad    
  - \int_{\Om} 
  \frac{ \gamma_{+} [ \Z^{\gamma_{+}-1} (\p_{R}\Z ) ](r_{\ep},q_{\ep})    }{\ep^2   r_{\ep}  } (R_{\ep}-R_{\ep}^{\#})
  (\vv-\vv^{\#}) \cdot \nabla R_{\ep}^{\#}   dx
    \\
  & \qquad    
  - \int_{\Om} 
  \frac{ \gamma_{+} [ \Z^{\gamma_{+}-1} (\p_{Q}\Z ) ](r_{\ep},q_{\ep})    }{\ep^2    q_{\ep} } (Q_{\ep}-Q_{\ep}^{\#}) 
   (\vv-\vv^{\#}) \cdot \nabla Q_{\ep}^{\#}   dx
    \\
  & \qquad  
  - \int_{\Om} 
  \frac{ \gamma_{+} [ \Z^{\gamma_{+}-1} (\p_{R}\Z ) ](r_{\ep},q_{\ep})    }{\ep^2   r_{\ep}  } (R_{\ep}-R_{\ep}^{\#})  
  (r_{\ep}-r_{\ep}^{\#}) \div \vu^{\#}    dx  \\
  & \qquad  
  - \int_{\Om} 
  \frac{ \gamma_{+} [ \Z^{\gamma_{+}-1} (\p_{Q}\Z ) ](r_{\ep},q_{\ep})    }{\ep^2    q_{\ep} } (Q_{\ep}-Q_{\ep}^{\#})  (q_{\ep}-q_{\ep}^{\#}) \div \vu^{\#}    dx  \\
   & \qquad  
  + \int_{\Om} 
  \ep^{-2} (R_{\ep}-R_{\ep}^{\#})  (\vu-\vu^{\#})
  \cdot 
  [\p_{RR} p(\Z)  \nabla r_{\ep}  +  \p_{RQ} p(\Z) \nabla q_{\ep}] (r_{\ep},q_{\ep}) dx  \\
   & \qquad  
  + \int_{\Om} 
  \ep^{-2} (Q_{\ep}-Q_{\ep}^{\#})  (\vu-\vu^{\#})
  \cdot 
  [\p_{QR} p(\Z)  \nabla r_{\ep}  +  \p_{QQ} p(\Z) \nabla q_{\ep}] (r_{\ep},q_{\ep}) dx \\
  & 
  \quad 
  \leq 
  \f{\mu}{100} \n{ \Grad (\vu-\vu^{\#})  }_{L^2}^2 
+ \f{\mu+\lambda}{100} \n{ \div (\vu-\vu^{\#})  }_{L^2}^2 \\
   & \qquad  +
   C \sp{
   \n{\ep^{-1}(R_{\ep}-R_{\ep}^{\#}) }_{L^2}^2
+
\n{\ep^{-1}(Q_{\ep}-Q_{\ep}^{\#}) }_{L^2}^2
+
\n{ \vu-\vu^{\#}}_{L^2}^2
   } \\
   & \qquad  +
   C 
   \sp{
\n{\ep^{-1}(r_{\ep}-r_{\ep}^{\#}) }_{L^2}^2
+
\n{\ep^{-1}(q_{\ep}-q_{\ep}^{\#}) }_{L^2}^2
+
\n{ \vv-\vv^{\#}}_{L^2}^2
}. 
\end{align}   }}
Applying Gr\"{o}nwall's inequality to the above inequality gives
\begin{align}
& \sup_{0 \leq t \leq T_{\ast}}
\sp{
\n{\ep^{-1}(R_{\ep}-R_{\ep}^{\#}) }_{L^2}^2
+
\n{\ep^{-1}(Q_{\ep}-Q_{\ep}^{\#}) }_{L^2}^2
+
\n{ \vu-\vu^{\#}}_{L^2}^2
} \\
& \quad 
\leq C \exp(C T_{\ast}) 
\int_0^t  \sp{
\n{\ep^{-1}(r_{\ep}-r_{\ep}^{\#}) }_{L^2}^2
+
\n{\ep^{-1}(q_{\ep}-q_{\ep}^{\#}) }_{L^2}^2
+
\n{ \vv-\vv^{\#}}_{L^2}^2
} d\tau \\
& \quad \leq 
C T_{\ast} \exp(C T_{\ast}) 
\sup_{0 \leq t \leq T_{\ast}} 
\sp{
\n{\ep^{-1}(r_{\ep}-r_{\ep}^{\#}) }_{L^2}^2
+
\n{\ep^{-1}(q_{\ep}-q_{\ep}^{\#}) }_{L^2}^2
+
\n{ \vv-\vv^{\#}}_{L^2}^2
} , 
\end{align}
which readily implies {\small{
\begin{align}
\int_0^t \n{  \vu-\vu^{\#} }_{H^1}^2 d\tau  
&
\leq C 
[ T_{\ast}^2 \exp(C T_{\ast})  +T_{\ast} ]  \\
& \times
\sup_{0 \leq t \leq T_{\ast}} 
\sp{
\n{\ep^{-1}(r_{\ep}-r_{\ep}^{\#}) }_{L^2}^2
+
\n{\ep^{-1}(q_{\ep}-q_{\ep}^{\#}) }_{L^2}^2
+
\n{ \vv-\vv^{\#}}_{L^2}^2
} . 
\end{align}   }}
This verifies \eqref{eq:cont-map} as long as we take $T_{\ast}$ suitably small such that
\begin{align}
 C 
[ T_{\ast}^2 \exp(C T_{\ast})  +T_{\ast} + T_{\ast} \exp(C T_{\ast}) ] 
<1.
\end{align}

\subsection{The end of proof of Theorem \ref{Thm1}}
The proof of Theorem \ref{Thm1} contains two parts: the first one concerns the uniform regularity of solutions, and the other one treats the convergence to the incompressible Navier-Stokes flows.

Let us start with the uniform regularity of solutions. To this end, we define 
\begin{align}
R_{\ep}^{(0)}:= R_{0,\varepsilon},\quad     
Q_{\ep}^{(0)}:= Q_{0,\varepsilon},\quad
\vu^{(0)}:= \mathbf{u}_{0,\varepsilon}.
\end{align}
Given $R_{\ep}^{(k)},Q_{\ep}^{(k)},\vu^{(k)}$ with $k \geq 0$, we define $R_{\ep}^{(k+1)},Q_{\ep}^{(k+1)},\vu^{(k+1)}$ through the following linear system:
\begin{align}\label{eq:appro}
    \begin{dcases}
        \partial_t R_{\ep}^{(k+1)} +\vu^{(k)} \cdot \Grad R_{\ep}^{(k+1)}+ R_{\ep}^{(k)} \div \vu^{(k+1)}= 0, \\
        \partial_t Q_{\ep}^{(k+1)} +\vu^{(k)} \cdot \Grad Q_{\ep}^{(k+1)}+ Q_{\ep}^{(k)} \div \vu^{(k+1)}= 0, \\
        \p_t \vu^{(k+1)}+ \vu^{(k)} \cdot \Grad \vu^{(k+1)} + 
        \frac{ \gamma_{+} [ \Z^{\gamma_{+}-1} (\p_{R}\Z ) ](R_{\ep}^{(k)},Q_{\ep}^{(k)})    }{\ep^2   (R_{\ep}^{(k)} + Q_{\ep}^{(k)}) }
        \Grad R_{\ep}^{(k+1)}  \\
          \qquad \qquad 
        +
        \frac{ \gamma_{+} [ \Z^{\gamma_{+}-1} (\p_{Q}\Z ) ](R_{\ep}^{(k)},Q_{\ep}^{(k)})    }{\ep^2   (R_{\ep}^{(k)} + Q_{\ep}^{(k)}) }
        \Grad Q_{\ep}^{(k+1)}  
       \\
        \qquad \qquad 
        = \f{\mu}{  R_{\ep}^{(k)}+Q_{\ep}^{(k)} } \Delta \vu^{(k+1)}      
        +\f{  \mu+\lambda  }{  R_{\ep}^{(k)}+Q_{\ep}^{(k)} } \nabla \div \vu^{(k+1)}, 
    \end{dcases}
\end{align}  
with initial data
\begin{align}
R_{\ep}(0,\cdot)=  R_{0,\varepsilon}, \quad 
Q_{\ep}(0,\cdot)=  Q_{0,\varepsilon}, \quad 
\mathbf{u}_{\ep}(0,\cdot)= \mathbf{u}_{0,\varepsilon}. 
\end{align}
By the same arguments as subsections \ref{subs:unif-est}-\ref{subs:clo}, we see for all $k \geq 0$ that 
\begin{align}
 R_{\ep}^{(k)},Q_{\ep}^{(k)},\vu^{(k)} \in L^{\infty}(0,T_{\ast}; H^s(\Om)) \cap \text{Lip} ([0,T_{\ast}]; H^{s-1}(\Om))   
\end{align}
with $s \geq 4$ and moreover the following uniform bounds hold:
\begin{align}\label{eq:unif-bdds}
    \begin{dcases}
        \mathscr{E}_{s} ( R_{\ep}^{(k)},Q_{\ep}^{(k)},\vu^{(k)})(t) 
        +
        \int_0^t \sp{
        \mu \n{ \nabla \vu^{(k)} }_{H^s}^2
        + (\mu+\lambda) \n{ \div \vu^{(k)} }_{H^s}^2
        } d\tau \leq C
        ,
       \\
        \mathscr{E}_{s-1} (\p_t R_{\ep}^{(k)}, \p_t Q_{\ep}^{(k)},\p_t \vu^{(k)})(t) \\
        \qquad \qquad 
        +
        \int_0^t \sp{
        \mu \n{ \nabla \p_t \vu^{(k)}}_{H^{s-1}}^2
        + (\mu+\lambda) \n{ \div \p_t \vu^{(k)} }_{H^{s-1}}^2
        } d\tau \leq C, 
    \end{dcases}
\end{align}
for all $t\in [0,T_{\ast}]$, provided that we take $\ep>0$ suitably small. 
Now, based on the contraction of the map established in subsection \ref{subs:cont}, it is a routine matter to pass to the limit $k\rightarrow \infty$ for the sequence $\{ (R_{\ep}^{(k)},Q_{\ep}^{(k)},\vu^{(k)}) \}_{k \geq 0}$ in order to identify the limit 
\begin{align}
  & R_{\ep}^{(k)} \rightarrow R_{\ep},   \quad
  Q_{\ep}^{(k)} \rightarrow Q_{\ep}  \quad \text{in}  \quad 
  L^{\infty}(0,T_{\ast}; L^2(\Om))\quad \text{as} \quad k\rightarrow \infty, \\
  & \vu^{(k)} \rightarrow  \vu  \quad \text{in}  \quad 
  L^{\infty}(0,T_{\ast}; L^2(\Om)) \cap L^2(0,T_{\ast};H^1(\Om))
  \quad \text{as} \quad k\rightarrow \infty.
\end{align} 
Then, it follows from \eqref{eq:unif-bdds} that $(R_{\ep},Q_{\ep},\vu)$ satisfies the uniform estimates \eqref{eq:unif-reg}. Clearly, $(R_{\ep},Q_{\ep},\vu)$ is a classical solution to \eqref{eq:two-fluid-1} with initial data $(R_{0,\ep},Q_{0,\ep},\mathbf{u}_{0,\ep})$ on the time interval $[0,T_{\ast}]$. Uniqueness of solutions is proved similar to \eqref{eq:cont-map}. The details are omitted.

Now we show that $(R_{\ep},Q_{\ep},\vu)$ converges to the local-in-time classical solution of the incompressible Navier-Stokes flows \eqref{eq:limit}, as $\ep \rightarrow 0$. It follows from the uniform estimates \eqref{eq:unif-reg} that up to a suitable subsequence 
 \begin{align}
    \begin{dcases}
        R_{\ep} \rightarrow 1 \quad \text{ strongly in } 
        L^{\infty}(0,T_{\ast}; H^s(\Om)) \cap \text{Lip} ([0,T_{\ast}]; H^{s-1}(\Om))
         ,
       \\
         Q_{\ep} \rightarrow 1 \quad \text{ strongly in } 
        L^{\infty}(0,T_{\ast}; H^s(\Om)) \cap \text{Lip} ([0,T_{\ast}]; H^{s-1}(\Om))
         .
       \\
        \vu \rightarrow \mathbf{u} \quad \text{ weakly-$\ast$ in } 
        L^{\infty}(0,T_{\ast}; H^s(\Om)) \cap \text{Lip} ([0,T_{\ast}]; H^{s-1}(\Om))
    \end{dcases}
\end{align}
for some $\mathbf{u} \in L^{\infty}(0,T_{\ast}; H^s(\Om)) \cap \text{Lip} ([0,T_{\ast}]; H^{s-1}(\Om))$. 
We further deduce from the well-known Aubin-Lions lemma that 
\begin{align}
     \vu \rightarrow \mathbf{u} \quad \text{ strongly in } 
        C ( [0,T_{\ast}]; H^{s-\delta'}(\Om) )
\end{align}
for some $\delta'>0$ sufficiently small. With the above convergence properties, it is easily verified that $\mathbf{u}$ solves the incompressible Navier-Stokes flows \eqref{eq:limit} with initial value $\mathbf{u}_0$.

\section{Proof of Theorem \ref{Thm2}}\label{Sec:4}
\subsection{Relative energy inequality} 
In this subsection, we establish the relative energy inequality by choosing the target solution as a test function. 
To this end, let $(R_{\ep},Q_{\ep},\vu)$ be the local solution to the primitive system \eqref{eq:two-fluid-1} and $\mathbf{u}$ be the local solution to the target system \eqref{eq:limit} on $[0,T_{\ast}]$, ensured by Theorem \ref{Thm1}.  
Following the computations in Bresch et al. \cite{Bre-Muc-Zat-19}, we obtain the following energy inequality {\small{
\begin{align}\label{eq:ener-ineq}
    & \int_{\Om}  \left[ 
    \f{1}{2} (R_{\ep}+Q_{\ep}) |\vu|^2   
   + \f{1}{\ep^2}
    \f{1}{\gamma_{+}-1} \sp{ \f{R_{\varepsilon}}{ \alpha_{\varepsilon} }}^{\gamma_{+}  } 
    \alpha_{\varepsilon} 
    +  \f{1}{\ep^2}
    \f{1}{\gamma_{-}-1} \sp{ \f{Q_{\varepsilon}}{ 1-\alpha_{\varepsilon} }}^{\gamma_{-}  } 
    \sp{ 1-\alpha_{\varepsilon} } 
    \right]
    dx \\
    & \qquad 
    + 
    \int_0^t \int_{\Om} 
    \Big(
    \mu |\nabla \vu|^2 +
    (\mu+\lambda) (\div \vu)^2
    \Big) 
    dx d\tau \\
    & \quad 
    \leq 
    \int_{\Om}  \left[ 
    \f{1}{2} (R_{0,\ep}+Q_{0,\ep}) |\mathbf{u}_{0,\ep}|^2   
   + \f{1}{\ep^2}
    \f{1}{\gamma_{+}-1} \sp{ \f{R_{0,\varepsilon}}{ \alpha_{0,\varepsilon} }}^{\gamma_{+}  } 
    \alpha_{0,\varepsilon} 
    +  \f{1}{\ep^2}
    \f{1}{\gamma_{-}-1} \sp{ \f{Q_{0,\varepsilon}}{ 1-\alpha_{0,\varepsilon} }}^{\gamma_{-}  } 
    \sp{ 1-\alpha_{0,\varepsilon} } 
    \right]
    dx. 
\end{align} }}
Here and in what follows, we denote by $\alpha_{\ep}:=R_{\ep}/Z_{\ep}$ and $\alpha_{0,\ep}:=R_{0,\ep}/Z_{0,\ep}$.  
For the target system \eqref{eq:limit}, it follows readily the energy inequality
\begin{align}\label{eq:ener-ineq-target}
    \int_{\Om} |\mathbf{u}|^2 dx 
    + 
    \mu \int_0^t \int_{\Om} | \nabla \mathbf{u}|^2 dx d\tau
     \leq 
  \int_{\Om} |\mathbf{u}_{0}|^2 dx . 
\end{align}
To proceed, we test the momentum equation \eqref{eq:two-fluid}$_3$ by $\mathbf{u}$ and obtain 
\begin{align}
& \int_{\Om} (R_{\ep}+Q_{\ep}) \vu \cdot \mathbf{u} \,dx
+ 
\int_0^t \int_{\Om} 
(R_{\ep}+Q_{\ep}) \vu \cdot 
\Big( 
\mathbf{u} \cdot \nabla \mathbf{u}+ \nabla \Pi-\f{\mu}{2} \Delta \mathbf{u}
\Big) \,dx d\tau \\
& \qquad
- \int_0^t \int_{\Om}
(R_{\ep}+Q_{\ep}) \vu \otimes \vu : \nabla \mathbf{u} \,dx d\tau
+ 
\mu 
\int_0^t \int_{\Om} \nabla \vu : \nabla \mathbf{u} \,dx d\tau \\
& \quad
=
\int_{\Om}
(R_{0,\ep}+Q_{0,\ep}) \mathbf{u}_{0,\ep} \cdot \mathbf{u}_0 \,dx , 
\end{align}
which implies via integration by parts that {\small{
\begin{align}\label{eq:momen}
&
  \int_{\Om} \sqrt{ R_{\ep}+Q_{\ep} }  \,\sqrt{2}  \, 
  \vu \cdot \mathbf{u} \,dx
  + 
  \int_{\Om}   \sqrt{ R_{\ep}+Q_{\ep} }   
  \sp{
  \sqrt{ R_{\ep}+Q_{\ep} }   - \sqrt{2}
  }
  \, 
  \vu \cdot \mathbf{u} \,dx \\
& \qquad  
+
\int_0^t \int_{\Om} 
(R_{\ep}+Q_{\ep}) \vu \cdot 
\Big( 
\mathbf{u} \cdot \nabla \mathbf{u}+ \nabla \Pi 
\Big) \,dx d\tau 
+ 
\f{\mu}{2} \int_0^t \int_{\Om} 
\Big( 
2-(R_{\ep}+Q_{\ep})
\Big) \vu \cdot \Delta \mathbf{u} \,dx d\tau 
\\
& \qquad 
+ 
2 \mu 
\int_0^t \int_{\Om} \nabla \vu : \nabla \mathbf{u} \,dx d\tau
- 
\int_0^t \int_{\Om}
(R_{\ep}+Q_{\ep}) \vu \otimes \vu : \nabla \mathbf{u} \,dx d\tau \\
& \quad 
= 
\int_{\Om} \sqrt{ R_{0,\ep}+Q_{0,\ep} }  \,\sqrt{2}  \, 
  \mathbf{u}_{0,\ep} \cdot \mathbf{u}_{0} \,dx
  + 
  \int_{\Om}   \sqrt{ R_{0,\ep}+Q_{0,\ep} }   
  \sp{
  \sqrt{ R_{0,\ep}+Q_{0,\ep} }   - \sqrt{2}
  }
  \, 
  \mathbf{u}_{0,\ep} \cdot \mathbf{u}_{0} \,dx . 
\end{align} }}
Combining \eqref{eq:ener-ineq}, \eqref{eq:ener-ineq-target} and \eqref{eq:momen}, we arrive at {\small{
\begin{align}\label{rel-ene-ine}
& \f{1}{2} 
\int_{\Om} \left|
\sqrt{ R_{\ep}+Q_{\ep} } \, \vu-  \sqrt{2} \,\mathbf{u}
\right|^2 dx \\
& \qquad 
+ 
\int_{\Om}  \left[  
    \f{1}{\ep^2}
    \f{1}{\gamma_{+}-1} \sp{ \f{R_{\varepsilon}}{ \alpha_{\varepsilon} }}^{\gamma_{+}  } 
    \alpha_{\varepsilon} 
    +  \f{1}{\ep^2}
    \f{1}{\gamma_{-}-1} \sp{ \f{Q_{\varepsilon}}{ 1-\alpha_{\varepsilon} }}^{\gamma_{-}  } 
    \sp{ 1-\alpha_{\varepsilon} } 
    \right]
    dx \\
    &  \qquad 
    + \mu \int_0^t \int_{\Om} \left| \nabla \vu-\nabla \mathbf{u} \right|^2  \, dx d\tau
    + (\mu+\lambda) \int_0^t \int_{\Om}
   (\div \vu)^2  \, dx d\tau \\
& \quad 
\leq 
\f{1}{2} 
\int_{\Om} \left|
\sqrt{ R_{0,\ep}+Q_{0,\ep} } \, \mathbf{u}_{0,\ep}-  \sqrt{2} \,\mathbf{u}_{0}
\right|^2 dx \\
& \qquad 
+ 
\int_{\Om}  \left[  
    \f{1}{\ep^2}
    \f{1}{\gamma_{+}-1} \sp{ \f{R_{0,\varepsilon}}{ \alpha_{0,\varepsilon} }}^{\gamma_{+}  } 
    \alpha_{0,\varepsilon} 
    +  \f{1}{\ep^2}
    \f{1}{\gamma_{-}-1} \sp{ \f{Q_{0,\varepsilon}}{ 1-\alpha_{0,\varepsilon} }}^{\gamma_{-}  } 
    \sp{ 1-\alpha_{0,\varepsilon} } 
    \right]
    dx \\
&\qquad 
    + 
\f{\mu}{2} \int_0^t \int_{\Om} 
\Big( 
2-(R_{\ep}+Q_{\ep})
\Big) \vu \cdot \Delta \mathbf{u} \,dx d\tau  \\
& \qquad 
+
\int_0^t \int_{\Om} 
(R_{\ep}+Q_{\ep}) \vu \cdot 
\Big( 
\mathbf{u} \cdot \nabla \mathbf{u}+ \nabla \Pi 
\Big) \,dx d\tau 
- 
\int_0^t \int_{\Om}
(R_{\ep}+Q_{\ep}) \vu \otimes \vu : \nabla \mathbf{u} \,dx d\tau \\
& \qquad 
+ 
\int_{\Om}   \sqrt{ R_{\ep}+Q_{\ep} }   
  \sp{
  \sqrt{ R_{\ep}+Q_{\ep} }   - \sqrt{2}
  }
  \, 
  \vu \cdot \mathbf{u} \,dx  \\
  & \qquad
  - 
  \int_{\Om}   \sqrt{ R_{0,\ep}+Q_{0,\ep} }   
  \sp{
  \sqrt{ R_{0,\ep}+Q_{0,\ep} }   - \sqrt{2}
  }
  \, 
  \mathbf{u}_{0,\ep} \cdot \mathbf{u}_{0} \,dx .
\end{align}  }}
Hence, it remains to estimate the right-hand side of \eqref{rel-ene-ine} suitably.

\subsection{Estimates for the remainder}
For the initial data, we see from the assumptions of Theorem \ref{Thm2} that 
\begin{align}\label{rem:init}
& 
 \f{1}{2} 
\int_{\Om} \left|
\sqrt{ R_{0,\ep}+Q_{0,\ep} } \, \mathbf{u}_{0,\ep}-  \sqrt{2} \,\mathbf{u}_{0}
\right|^2 dx   \\
& \quad 
\leq  
C \int_{\Om} \left|
\sqrt{ R_{0,\ep}+Q_{0,\ep} } \, \mathbf{u}_{0,\ep}-  \sqrt{2} \,\mathbf{u}_{0,\ep}
\right|^2 dx
+
C \int_{\Om} \left|
\mathbf{u}_{0,\ep} -\mathbf{u}_{0}
\right|^2 dx \\
& \quad
\leq 
C \int_{\Om}
\f{  \left| R_{0,\ep}-1 \right|^2 +\left| Q_{0,\ep}-1 \right|^2     }{  \left| \sqrt{ R_{0,\ep}+Q_{0,\ep} } +  \sqrt{2} \right|^2    }
| \mathbf{u}_{0,\ep}|^2 dx
+ 
C \int_{\Om} \left|
\mathbf{u}_{0,\ep} -\mathbf{u}_{0}
\right|^2 dx \\
& \quad
\leq C \ep ,
 \\ 
& 
   \int_{\Om}  \left[  
    \f{1}{\ep^2}
    \f{1}{\gamma_{+}-1} \sp{ \f{R_{0,\varepsilon}}{ \alpha_{0,\varepsilon} }}^{\gamma_{+}  } 
    \alpha_{0,\varepsilon} 
    +  \f{1}{\ep^2}
    \f{1}{\gamma_{-}-1} \sp{ \f{Q_{0,\varepsilon}}{ 1-\alpha_{0,\varepsilon} }}^{\gamma_{-}  } 
    \sp{ 1-\alpha_{0,\varepsilon} } 
    \right]
    dx 
    \leq 
    C \ep. 
\end{align}
Next, by H\"{o}lder's inequality and the uniform estimates \eqref{eq:unif-reg}, it holds
\begin{align}\label{rem:den}
&
\f{\mu}{2} \int_0^t \int_{\Om} 
\Big( 
2-(R_{\ep}+Q_{\ep})
\Big) \vu \cdot \Delta \mathbf{u} \,dx d\tau \\
& \quad 
\leq 
C 
\sp{
\n{R_{\ep}-1}_{L^2(0,t;L^2)}
+ \n{Q_{\ep}-1}_{L^2(0,t;L^2)}
}
\n{ \vu}_{L^2(0,t;L^2)}
\n{\Delta \mathbf{u}}_{L^{\infty}(0,t;L^{\infty})} \\
& \quad 
\leq C \ep. 
\end{align}
By the continuity equations \eqref{eq:two-fluid}$_1$-\eqref{eq:two-fluid}$_2$ and integration by parts, 
\begin{align}\label{rem:dens}
 \int_0^t \int_{\Om} 
 \sp{ R_{\ep}+Q_{\ep} } \vu \cdot \nabla \Pi \, dxd\tau
 & =
  \int_0^t \int_{\Om} 
  \Pi \, \p_t  \sp{ R_{\ep}+Q_{\ep} }    \, dxd\tau
  \\
  & \leq 
  \n{   \p_t  \sp{ R_{\ep}+Q_{\ep} }  }_{L^2(0,t;L^2)}
  \n{\Pi}_{L^2(0,t;L^2)}
  \leq C \ep. 
\end{align}
To continue, we have 
\begin{align}
 &   \int_{\Om}   \sqrt{ R_{\ep}+Q_{\ep} }   
  \sp{
  \sqrt{ R_{\ep}+Q_{\ep} }   - \sqrt{2}
  }
  \, 
  \vu \cdot \mathbf{u} \,dx  \\
& \quad 
\leq 
C 
\n{ \sqrt{ R_{\ep}+Q_{\ep} }   - \sqrt{2} }_{L^2}
\n{  \sqrt{ R_{\ep}+Q_{\ep} } \,\vu  }_{L^2}
\n{ \mathbf{u}}_{L^{\infty}} \\
& \quad 
\leq 
C  \sp{ 
\n{R_{\ep}-1 }_{L^2} + \n{Q_{\ep}-1 }_{L^2}
}
\n{  \sqrt{ R_{\ep}+Q_{\ep} } \,\vu  }_{L^2}
\n{ \mathbf{u}}_{L^{\infty}} \\
& \quad 
\leq C \ep.
\end{align}
A similar estimate holds for the integral involved with the initial data.

The remaining two integrals require a bit more effort. We first rearrange as 
\begin{align}
 & \int_0^t \int_{\Om} 
(R_{\ep}+Q_{\ep}) \vu \cdot 
\Big( 
\mathbf{u} \cdot \nabla \mathbf{u}
\Big) \,dx d\tau 
- 
\int_0^t \int_{\Om}
(R_{\ep}+Q_{\ep}) \vu \otimes \vu : \nabla \mathbf{u} \,dx d\tau  \\
& \quad =
- \int_0^t \int_{\Om}
\sp{ 
\sqrt{ R_{\ep}+Q_{\ep}} \, \vu - \sqrt{2}\, \mathbf{u}
} 
\otimes 
\sp{ 
\sqrt{ R_{\ep}+Q_{\ep}} \, \vu - \sqrt{2}\, \mathbf{u}
}: \nabla \mathbf{u}  \,dx d\tau  \\
& \qquad
+
\int_0^t \int_{\Om}
\sp{
R_{\ep}+Q_{\ep}- \sqrt{ R_{\ep}+Q_{\ep}} \, \sqrt{2}
} \vu \cdot \Big( 
\mathbf{u} \cdot \nabla \mathbf{u}
\Big) \,dx d\tau \\
& \qquad 
-
\int_0^t \int_{\Om} 
\sp{
\sqrt{ R_{\ep}+Q_{\ep}} \, \sqrt{2} \, \vu- 2 \mathbf{u}
} \cdot \nabla \sp{ \f{|\mathbf{u}|^2}{2}   } 
 \,dx d\tau, 
\end{align}
while for the last term we further rewrite, by the divergence-free condition \eqref{eq:limit}$_2$ and the continuity equations \eqref{eq:two-fluid}$_1$-\eqref{eq:two-fluid}$_2$, as {\small{
\begin{align}
&
\int_0^t \int_{\Om} 
\sp{
\sqrt{ R_{\ep}+Q_{\ep}} \, \sqrt{2} \, \vu- 2 \mathbf{u}
} \cdot \nabla \sp{ \f{|\mathbf{u}|^2}{2}   } 
 \,dx d\tau 
 \\
 & \quad 
 = 
 \int_0^t \int_{\Om} 
\sp{
\sqrt{ R_{\ep}+Q_{\ep}} \, \sqrt{2} \, \vu
-( R_{\ep}+Q_{\ep}) \vu+ ( R_{\ep}+Q_{\ep}) \vu
- 2 \mathbf{u}
} \cdot \nabla \sp{ \f{|\mathbf{u}|^2}{2}   } 
 \,dx d\tau \\
  & \quad 
 = 
 \int_0^t \int_{\Om} 
 \sqrt{ R_{\ep}+Q_{\ep}}
 \sp{
 \sqrt{2}- \sqrt{ R_{\ep}+Q_{\ep}}
 } \vu \cdot \sp{ \f{|\mathbf{u}|^2}{2}   } 
 \,dx d\tau \\
 & \qquad \qquad 
 +
  \int_0^t \int_{\Om} 
   \f{|\mathbf{u}|^2}{2}
   \p_t (R_{\ep}+Q_{\ep}) \,dx d\tau.
\end{align}  }}
From the two equalities above and the uniform estimates \eqref{eq:unif-reg}, we readily get 
\begin{align}\label{rem:conv}
&
  \left|
\int_0^t \int_{\Om} 
(R_{\ep}+Q_{\ep}) \vu \cdot 
\Big( 
\mathbf{u} \cdot \nabla \mathbf{u}
\Big) \,dx d\tau 
- 
\int_0^t \int_{\Om}
(R_{\ep}+Q_{\ep}) \vu \otimes \vu : \nabla \mathbf{u} \,dx d\tau 
  \right|  \\
  & \quad \leq 
  C \int_0^t \int_{\Om} 
\left|
\sqrt{ R_{\ep}+Q_{\ep} } \, \vu-  \sqrt{2} \,\mathbf{u}
\right|^2 dx d\tau 
+ C \ep. 
\end{align}

Based on \eqref{rem:init}, \eqref{rem:den}, \eqref{rem:dens} and \eqref{rem:conv}, we see from \eqref{rel-ene-ine} that 
\begin{align}
& \f{1}{2} 
\int_{\Om} \left|
\sqrt{ R_{\ep}+Q_{\ep} } \, \vu-  \sqrt{2} \,\mathbf{u}
\right|^2 \,dx \\
& \qquad 
+ 
\int_{\Om}  \left[  
    \f{1}{\ep^2}
    \f{1}{\gamma_{+}-1} \sp{ \f{R_{\varepsilon}}{ \alpha_{\varepsilon} }}^{\gamma_{+}  } 
    \alpha_{\varepsilon} 
    +  \f{1}{\ep^2}
    \f{1}{\gamma_{-}-1} \sp{ \f{Q_{\varepsilon}}{ 1-\alpha_{\varepsilon} }}^{\gamma_{-}  } 
    \sp{ 1-\alpha_{\varepsilon} } 
    \right]
    \,dx \\
    &  \qquad 
    + \mu \int_0^t \int_{\Om} \left| \nabla \vu-\nabla \mathbf{u} \right|^2  \, dx d\tau
    + (\mu+\lambda) \int_0^t \int_{\Om}
   (\div \vu)^2  \, dx d\tau \\
   & \quad 
\leq 
C \int_0^t \int_{\Om} 
\left|
\sqrt{ R_{\ep}+Q_{\ep} } \, \vu-  \sqrt{2} \,\mathbf{u}
\right|^2  \, dx d\tau 
+ C \ep, 
\end{align}
from which we infer by Gr\"{o}nwall's inequality that 
\begin{align}\label{eq:bdd-rel-ene}
& \f{1}{2} 
\int_{\Om} \left|
\sqrt{ R_{\ep}+Q_{\ep} } \, \vu-  \sqrt{2} \,\mathbf{u}
\right|^2 dx \\
& \qquad 
+ 
\int_{\Om}  \left[  
    \f{1}{\ep^2}
    \f{1}{\gamma_{+}-1} \sp{ \f{R_{\varepsilon}}{ \alpha_{\varepsilon} }}^{\gamma_{+}  } 
    \alpha_{\varepsilon} 
    +  \f{1}{\ep^2}
    \f{1}{\gamma_{-}-1} \sp{ \f{Q_{\varepsilon}}{ 1-\alpha_{\varepsilon} }}^{\gamma_{-}  } 
    \sp{ 1-\alpha_{\varepsilon} } 
    \right]
    dx \\
    &  \qquad 
    + \mu \int_0^t \int_{\Om} \left| \nabla \vu-\nabla \mathbf{u} \right|^2  \, dx d\tau
    + (\mu+\lambda) \int_0^t \int_{\Om}
   (\div \vu)^2  \, dx d\tau \\
   & \quad 
\leq 
 C \ep.
\end{align}

\subsection{Convergence rates}
To begin with, we observe that 
\begin{align}\label{rate:den}
    \n{R_{\ep}-1}_{H^s}^2 + \n{Q_{\ep}-1}_{H^s}^2
\leq C \ep^2
\end{align}
follows directly from the uniform estimates \eqref{eq:unif-reg}. Next, we see from \eqref{eq:bdd-rel-ene} that 
\begin{align}\label{rate:vel}
\n{
\vu-\mathbf{u}
}_{L^2}^2 
& 
=\f{1}{2}
\n{
\sqrt{2} \,\vu-  \sqrt{2} \, \mathbf{u}
}_{L^2}^2 \\
& 
\leq 
C \sp{
\n{
\sqrt{ R_{\ep}+Q_{\ep} }\, \vu 
- 
\sqrt{2} \,  \vu
}_{L^2}^2 
+ 
\n{
\sqrt{ R_{\ep}+Q_{\ep} }\, \vu 
- 
\sqrt{2} \,  \mathbf{u}
}_{L^2}^2 
} \\
& 
\leq 
C \ep^2 + C \ep \leq C \ep, 
\end{align}
where we used 
\begin{align}
   \n{
\sqrt{ R_{\ep}+Q_{\ep} }\, \vu 
- 
\sqrt{2} \,  \vu
}_{L^2}^2 
\leq 
C 
\n{\vu}_{L^{\infty} }^2
\sp{ 
\n{R_{\ep}-1}_{L^2}^2 
+
\n{Q_{\ep}-1}_{L^2}^2 
}
\leq 
C \ep^2. 
\end{align}
To proceed, we recall from \eqref{eq:two-fluid}$_1$ that 
\begin{align}
\div \vu=
- \sp{
\p_t R_{\ep} + \sp{ R_{\ep}-1} \div \vu 
+ \vu \cdot \Grad R_{\ep}
}; 
\end{align}
whence Lemma \ref{lemm:prod} and the uniform estimates \eqref{eq:unif-reg} yield 
\begin{align}
    \n{ \ep^{-1}  \div \vu }_{H^{s-1}}
   & \leq 
    C \Big(
    \n{  \ep^{-1} \p_t R_{\ep}  }_{H^{s-1}}
    + 
    \n{
   \ep^{-1}  \sp{ R_{\ep}-1} \div \vu
    }_{H^{s-1}}  \\
    & \qquad \qquad
    +
    \n{
    \ep^{-1}  \vu \cdot \Grad \sp{ R_{\ep} -1} 
    }_{H^{s-1}}
    \Big) \\
    & \leq 
    C 
    \sp{
    \n{  \ep^{-1} \p_t R_{\ep}  }_{H^{s-1}}
    + 
    \n{  \ep^{-1}  \sp{ R_{\ep}-1}  }_{H^{s}}
    } \\
    & \leq  C. 
\end{align}
This obviously shows that 
\begin{align}\label{rate:div-vel}
    \n{ \div \vu}_{H^{s-1}} 
= 
\n{ \div \sp{ \vu- \mathbf{u}}  }_{H^{s-1}} 
\leq C \ep. 
\end{align}
Putting \eqref{eq:bdd-rel-ene}, \eqref{rate:den}, \eqref{rate:vel} and \eqref{rate:div-vel} together, we obtain the convergence rates \eqref{eq:rate}, thus finishing the proof of Theorem \ref{Thm2}.

\vspace{10mm}

\noindent{{\bf{Acknowledgements}}} 
The work of Y.L. was supported by National Natural Science Foundation of China (12571228), Natural Science Foundation of Anhui Province (2408085MA018). The work of M.L.-M. was supported by the Gutenberg Research College and by the Deutsche Forschungsgemeinschaft (DFG, German Research Foundation)--project number 233630050--TRR 146 and project number 525853336--SPP 2410 ``Hyperbolic Balance Laws: Complexity, Scales and Randomness". She is also grateful to the Mainz Institute of Multiscale Modelling for supporting her research. The work of E.Z. was supported by the EPSRC Early Career Fellowship no. EP/V000586/1.

\vspace{4mm}

\noindent{{\bf{Data Availability}}} Data sharing is not applicable to this article as no datasets were generated or analysed
during the current study.

\vspace{4mm}

\noindent{{\bf{Conflicts of interest}}} All authors certify that there are no conflicts of interest for this work.

\vspace{4mm}

\noindent{\bf{Publishing licence.}} For the purpose of open access, the author has applied a Creative Commons Attribution (CC BY) licence to any Author Accepted Manuscript version arising from this submission.

\appendix

\def\thesection{\Alph{section}}\label{Sec:A}
\section{Two useful lemmas}


In this appendix, we recall some useful lemmas which are extensively applied in the proof of main theorem.

To begin with, we recall the following classical product estimate and the commutator estimate in Sobolev spaces from \cite{Kla-Maj-81}*{Lemma A.1}.
\begin{lemm}\label{lemm:prod}
Let $f,g \in H^{N}(\Om)$ with $N \geq 1$ being an integer. Then, for any multi-index $\alpha=(\alpha_1,\alpha_2,\alpha_3)$ with $|\alpha|\leq N$ we have
\begin{align}
& \n{ \nabla^{\alpha}(fg) }_{L^2} 
\leq C \sp{ 
\n{f}_{L^{\infty}}  \n{\nabla^{N} g}_{L^2} + \n{g}_{L^{\infty}} \n{\nabla^{N} f}_{L^2} 
}, \\
& \n{ \nabla^{\alpha}(fg)   -f \nabla^{\alpha} g }_{L^2}
\leq C \sp{
\n{ \nabla f }_{L^{\infty}} \n{ \nabla^{N-1} g }_{L^2}
+ \n{g}_{L^{\infty}} \n{\nabla^{N} f}_{L^2} 
},
\end{align}
for some positive constant $C$ depending only on $N$. 
\end{lemm}

Next, we recall the composition lemma in Sobolev spaces from \cite{Yao-Zhu-Zi-12}*{Lemma 2.2}. 
\begin{lemm}\label{lemm:compo}
Let $F:\R \rightarrow \R$ and $G:\R^2 \rightarrow \R$ be smooth functions. Assume that $f,g \in H^{N}(\Om) \cap L^{\infty}(\Om)$ with $N \geq 1$ being an integer. Then it holds 
\begin{align}
& \n{ \nabla F(f) }_{H^{N-1}}
\leq C \sp{
\n{\nabla f}_{H^{N-1}} +  \n{\nabla f}_{H^{N-1}}^N
}, \\
& \n{ \nabla G(f,g) }_{H^{N-1}}
\leq C
\sp{
\n{\nabla f}_{H^{N-1}} + \n{\nabla g}_{H^{N-1}}
+ \n{\nabla f}_{H^{N-1}}^N
+\n{\nabla g}_{H^{N-1}}^N
},
\end{align}
for some positive constant $C$ depending only on $N,F,G$ and $\n{(f,g)}_{L^{\infty}}$. 
\end{lemm}

\end{document}